\documentclass[a4paper,12pt,leqno]{article}
\usepackage{amsfonts}
\usepackage{amsmath}
\usepackage{amssymb}
\usepackage{amsthm}

\usepackage{color}

\usepackage[margin=20truemm]{geometry}

\usepackage[hyphens]{url}
\usepackage{hyperref}
\usepackage[hyphenbreaks]{breakurl}

\newcommand\blfootnote[1]{
  \begingroup
  \renewcommand\thefootnote{}\footnote{#1}
  \addtocounter{footnote}{-1}
  \endgroup
}

\usepackage{titlesec}
\titleformat{\section}{\centering\normalfont\normalfont\bfseries}{\thesection.}{2mm}{}

\usepackage{indentfirst}

\numberwithin{equation}{section}

\usepackage{hyperref}
\usepackage[noabbrev,capitalise,nameinlink]{cleveref}
\hypersetup{
    colorlinks=true,
    citecolor=blue,
    linkcolor=blue,
    urlcolor=blue,
}

\usepackage[
    uniquelist=false,
    style=alphabetic,
    maxnames=99,
    maxalphanames=4,
    minalphanames=3,
    abbreviate=true,
    sorting=nty,
    backref=true
]{biblatex}

\bibliography{refs}

\renewbibmacro{in:}{}

\DefineBibliographyStrings{english}{
  backrefpage  = {\hspace{-0.15cm}},
  backrefpages = {\hspace{-0.15cm}},
}

\newtheorem{thm}{Theorem}
\newtheorem{rem}{Remark}
\newtheorem{eg}{Example}

\newtheorem{THM}{Main Theorem}
\newenvironment{THM*}[1]{\renewcommand\theTHM{#1}\THM}{\endTHM}
\newtheorem{thmMain}{Theorem}
\newenvironment{thm*}[1]{\renewcommand\thethmMain{#1}\thmMain}{\endthmMain}
\newtheorem{propMain}{Proposition}
\newenvironment{prop*}[1]{\renewcommand\thepropMain{#1}\propMain}{\endpropMain}
\newtheorem{lemMain}{Lemma}
\newenvironment{lem*}[1]{\renewcommand\thelemMain{#1}\lemMain}{\endlemMain}

\newcommand{\Q}{\mathbb Q}
\newcommand{\Z}{\mathbb Z}

\def\GL{\mathrm{GL}}
\def\Ker{\mathrm{Ker}}

\def\:{\colon}

\begin{document}

\begin{center}
 {\Large \bf 
 On the Burde--de Rham theorem for \\ 
 finitely presented pro-$p$ groups 
 }
\end{center}

\begin{center}
Yasushi MIZUSAWA, 
Ryoto TANGE, and 
Yuji TERASHIMA
\end{center}



\noindent
{\bf Abstract.}
We consider the Burde--de Rham theorem for finitely presented pro-$p$ groups 
under the assumption that the total degrees of all relators are $0$. 
We also give some concrete examples including 
higher-dimensional 
cases under Iwasawa theoretic conditions, 
and consider some cohomological interpretations.

\section{Introduction}
In this paper, 
we study linear representations of finitely presented pro-$p$ groups from the knot theoretic viewpoint. 
In knot theory, 
Burde and de Rham independently discovered the constructions of linear representations of 
knot groups using the zeros of Alexander polynomials 
(\cite{burde1967Knoten, deRham1967Noeud, kitano2013linear}). 

The analogy between the Alexander polynomial and 
the algebraic $p$-adic $L$-function was 
pointed out by Mazur (\cite{Mazur1964}), 
and based on this analogy, 
various studies have developed. 
For these studies, see Chapters 10 – 15 of \cite{morishita2024knots} and 
the references therein. 
In knot theory, Alexander polynomials are extended to twisted Alexander polynomials 
for any representations of any dimension of the knot group and are actively studied
(\cite{lin2001representations, wada1994twisted}). 
Alexander polynomials are regarded as twisted Alexander polynomials
for the trivial one-dimensional representation. 
Twisted Alexander polynomials are defined based on 
the Fox free differential calculus (\cite{wada1994twisted}). 
From the viewpoint of arithmetic topology \cite{morishita2024knots}, 
we introduce arithmetic twisted Alexander polynomials 
based on the pro-$p$ Fox free differential calculus 
(\cite{ihara1986galois, morishita2002certainKP})
and show the arithmetic Burde--de Rham theorem 
which describes a relation between the zeros of the arithmetic twisted
Alexander polynomials and
extensions of representations for pro-$p$ Galois groups of number fields.

From now on, 
let us briefly mention 
the number theoretic analogue of the Burde--de Rham Theorem 
(\cite{burde1967Knoten, deRham1967Noeud}). 
Let $F_{n}$ be the free pro-$p$ group on $g_{1}, \dots, g_{n}$.  
The total degree of $r \in F_n$ is defined as the image by the homomorphism 
$\deg \: F_n \rightarrow \mathbb Z_p$ which sends all $g_i$ to $1$. 
Let $G$ be a finitely presented pro-$p$ group with 
$n$ generators and 
with the total degrees of all relators being $0$, 
Let $\pi \: F_{n} \to G$ be the natural homomorphism, and 
denote the same $g_{i}$ for the image of $g_{i}$ in $G$. 
Denote $\Gamma$ the ring of $p$-adic integers $\mathbb{Z}_{p}$, and 
suppose that we can take a surjective homomorphism 
$\alpha \: G \to \Gamma$ such that 
$\alpha(g_{i}) = \gamma$ 
for all $i$, 
where $\gamma$ is a topological generator of $\Gamma$. 
Let $\overline{\Q_{p}}$ be the algebraic closure of the field of $p$-adic numbers. 
For a given $p$-adic number $a_{1} \in \overline{\Q_{p}}^{\times}$ with 
the $p$-adic absolute value $|a_{1}-1|_p<1$, 
let $\overline{\beta_{1}} \: F_{n} \to \overline{\Q_{p}}$ 
be a crossed homomorphism with respect to 
$\widetilde{\varphi_{0}} |_{\gamma = a_{1}} \: F_{n} \to \overline{\Q_{p}}^{\times}$, 
where $\widetilde{\varphi_{0}} := 
\varphi_{0} \otimes (\alpha \circ \pi) \: F_{n} \to 
(\overline{\Q_p} \otimes_{\Z_p} \Z_p[[\Gamma]] )^{\times}$
is the tensor product representation 
of a given continuous homomorphism $\varphi_{0} \: F_{n} \to \overline{\Q_{p}}^{\times}$. 
For a given continuous homomorphism $\varphi_{1} \: F_{n} \to \GL(2; \overline{\Q_{p}} )$, 
we say $\varphi_{1}$ 
{\it factors through} 
$G$ as a representation 
$\rho_{1} \: G \to \GL(2; \overline{\Q_{p}} )$ 
if $\rho_{1} \circ \pi = \varphi_{1}$. 
Then a $2$-dimensional version of our main result is the following: 

\begin{thm*} {\hspace{-1.2mm}}
The continuous homomorphism defined by 
$$
\varphi_{1} \: F_{n} \ni g_{i} \mapsto 
\begin{pmatrix} 
\widetilde{ \varphi_{0} } (g_{i})|_{\gamma = a_{1}} & 
\overline{\beta_{1}}(g_{i}) \\ 
0 & 1 
\end{pmatrix} 
\in \GL(2; \overline{\Q_{p}} )
$$ 
for all $i$, can be 
factors through 
$G$ as a representation 
$\rho_{1} \: G \to \GL(2; \overline{\Q_{p}} )$ 
in at least two ways if and only if $a_{1} \in \overline{\Q_{p}}^{\times}$ is a zero of 
the arithmetic twisted Alexander polynomial, 
which is a characteristic polynomial of the Iwasawa module 
$(\mathrm{Ker} (\alpha) )^{{\rm ab}}$. 
\end{thm*}

\noindent
Note that $(\mathrm{Ker} (\alpha) )^{\rm ab}$ is the maximal abelian pro-$p$ quotient of 
$\mathrm{Ker} (\alpha)$, 
which is a module finitely generated over Iwasawa algebra $\Z_{p}[[\Gamma]]$. 

On the other hand in knot theory, 
generalizations of the Burde--de Rham Theorem to higher-dimensional linear representations are 
also studied 
(\cite{SilverWilliams2011theorem}), 
and are still actively studied from cohomological viewpoints 
(\cite{HeusenerPorti2015rep, HeusenerPortiSuarez2001defoRed}). 
In this paper, 
we will prove 
the number theoretic analogue of the Burde--de Rham Theorem 
including higher-dimensional cases (\cref{thm:main}), 
and give some cohomological interpretations (\cref{thm:zeroesCoh}, \cref{thm:cohZeroes}). 

One of the differences between knot groups and 
finitely presented pro-$p$ Galois groups of number fields is the deficiency of group presentations, 
namely the difference in the number of generators and relators of group presentations. 
The knot groups can always transform into having the deficiency one 
(called the Wirtinger presentation), whereas 
the pro-$p$ Galois groups cannot. 
In this paper, 
we consider the case of group presentations for any deficiency. 

Another concern is the number theoretic counterpart of 
the Alexander polynomial for higher-dimensional cases. 
Based on the ideas of twisted Alexander invariant in knot theory, 
we consider the counterpart using Fitting ideals. 
We also give some cohomological viewpoints, which we expect 
the relations with infinitesimal deformations of representations 
(\cite{mazur1989deforming}). 
We believe that this study will enrich the field of number theory 
(e.g., Iwasawa theory, the zeros of $p$-adic $L$-functions, etc.), 
and strengthen the bridge between number theory and knot theory 
(e.g., multiple residue symbols, Massey products, etc. 
(\cite{ahlqvistmagnus2022massey, morishita2004milnor, sharifi2007iwasawa, sharifi2007massey})). 

This paper is organized as follows. 
In \cref{sec:DiffMod}, 
we recall the complete differential modules and 
introduce some properties that we use later. 
In \cref{sec:mainThm}, 
we prove the Burde--de Rham theorem for finitely presented pro-$p$ groups 
including higher-dimensional cases, and 
consider the existence of non-standard extensions of homomorphisms 
under the assumption that the total degrees of all relators are $0$. 
In \cref{sec:Eg}, 
we give some explicit examples of extensions of homomorphisms for 
Galois groups under Iwasawa theoretic settings. 
Lastly, in \cref{sec:cohInterpret}, 
we recall the group cohomology of local systems and 
study some cohomological interpretations of extensions of homomorphisms, and 
consider generalizations of the relation with the zeros of algebraic $p$-adic $L$-functions. 
\cref{sec:cohInterpret} is seen as a sequel of the works by Morishita 
(\cite{morishita2006alexander}, \cite[Chapter 14]{morishita2024knots}), 
where the knot theoretic part is mainly due to 
Hironaka (\cite{hironaka1997alexander}) and L\^{e} (\cite{le1994varieties}). 

\ 

\noindent
\textbf{Notation. }
Throughout this paper, we only consider homomorphisms with continuity. 
For a pro-$p$ group $G$, 
we denote by $[h, g] = h^{-1} g^{-1} h g$ the commutator of $g, \ h \in G$. 
We denote by $|\cdot|_{p}$ the $p$-adic absolute value 
normalized by $|p|_{p} = p^{-1}$. 
For a closed subgroup $H$ of $G$, $[H, G]$ 
denotes the minimal closed subgroup containing $\{ [h, g] \mid h \in H, g \in G \}$. 
For $a$, $b$ in a commutative ring $R$, 
$a \doteq b$ means $a = b u$ for some unit 
$u \in R^{\times}$. 

\section{Complete differential modules} \label{sec:DiffMod}

Let $G$ and $H$ be profinite groups, and 
let $\psi \: G \to H$ be a continuous homomorphism.
Let $p$ be a prime number. 
For simplicity, 
we also denote by the same $\psi$ the algebra homomorphism 
$\mathbb{Z}_{p}[[G]] \to \mathbb{Z}_{p}[[H]]$ of the complete group algebras over 
the ring of $p$-adic integers $\mathbb{Z}_{p}$ induced by $\psi$. 

The {\it complete $\psi$-differential module} $A_{\psi}$ is defined as the quotient
module of the left free $\mathbb{Z}_{p}[[H]]$-module 
$\oplus_{g \in G} \mathbb{Z}_{p}[[H]]dg$ on the symbols
$dg$ ($g \in G$) by the left $\mathbb{Z}_{p}[[H]]$-submodule generated by elements of the form
$d(x y) - dx - \psi(x) dy$ for any $x$, $y \in G$:
$$
A_{\psi} :=  
( \oplus_{g \in G} \mathbb{Z}_{p}[[H]]dg ) / 
\langle
d(x y) - dx - \psi(x) dy \ (x, y \in G)
\rangle_{\mathbb{Z}_{p}[[H]]}.
$$

\noindent
The map $d \: G \to A_{\psi}$ defined by 
$g \mapsto dg$ is a $\psi$-differential, namely, 
we have 
$$d(x y) = dx + \psi(x) dy$$ 
for any $x$, $y \in G$, 
and the following universal property holds:

\ 

For any left $\mathbb{Z}_{p}[[H]]$-module $M$, and 
for any $\psi$-differential $\partial \: G \to M$, 

there exists a unique $\mathbb{Z}_{p}[[H]]$-homomorphism $f \: A_{\psi} \to M$ 
such that $f \circ d = \partial$. 

\ 

\noindent
Next, suppose that $G$ is a finitely presented pro-$p$ group, and 
choose a presentation 
\begin{equation} \label{eq:grpPres}
G = 
\langle g_{1}, \dots, g_{n} \mid 
r_1 = \cdots = r_m = 1 \rangle^{{\rm pro-}p}.     
\end{equation}
Let $F_{n}$ be the free pro-$p$ group on $g_{1}, \dots, g_{n}$, and 
let $\pi \: F_{n} \to G$ be the natural homomorphism. 
We also denote by 
the same $g_{i}$ for the image of $g_{i}$ in $G$, and 
the same $\pi$ for the algebra homomorphism 
$\mathbb{Z}_{p}[[F_{n}]] \to \mathbb{Z}_{p}[[G]]$ of complete group algebras 
induced by $\pi$. 
Then we can describe the presentation matrix of 
the complete $\psi$-differential module $A_{\psi}$ using 
the pro-$p$ Fox free differential calculus 
$\frac{\partial}{ \partial g_{i} } \: \mathbb{Z}_{p}[[F_{n}]] \to \mathbb{Z}_{p}[[F_{n}]]$ 
of complete group algebras 
(\cite[Section 2]{ihara1986galois}, 
\cite[Section 2]{morishita2002certainKP}, 
\cite[Section 10.3]{morishita2024knots}, 
\cite[Section V.6]{neukirch2020cohomology}), as follows:

\begin{prop*}{\hspace{-1.2mm}}[{\cite[Corollary 10.3.6]{morishita2024knots}}] \label{prop*:presMat} 
The complete $\psi$-differential module $A_{\psi}$ has a free resolution over $\mathbb{Z}_{p}[[H]]$:
$$
\mathbb{Z}_{p}[[H]]^{m} 
\stackrel{Q_{\psi \circ \pi}}{\longrightarrow} 
\mathbb{Z}_{p}[[H]]^{n} \longrightarrow A_{\psi} \to 0 
$$
whose presentation matrix 
$Q_{\psi \circ \pi} \in 
{\rm M}(m, n; \mathbb{Z}_{p} [[H]] )$ is given by
$$
Q_{\psi \circ \pi} := 
\left( (\psi \circ \pi) 
\left( \frac{ \partial r_{j} }{ \partial g_{i} } \right) 
\right)
\in {\rm M}(m, n; \mathbb{Z}_{p} [[H]] ). 
$$ 
\end{prop*}

For later use, we remark on some basic rules for Fox derivatives (\cite[Section 2.1]{ihara1986galois}): 
\begin{enumerate}
    \item 
    $\displaystyle
    \frac{\partial g_j}{\partial g_i} = 
    \left\{\begin{array}{ll} 1 & (i = j) \\ 0 & (i \neq j) \end{array}\right. ; 
    $ \\
    \item 
    $\displaystyle
    \frac{\partial (\eta + \xi)}{\partial g_i} = 
    \frac{\partial \eta}{\partial g_i} + \frac{\partial \xi}{\partial g_i}, \  
    \quad (\eta, \ \xi \in \Z_{p}[[F_{n}]]); 
    $ \\ 
    \item
    $\displaystyle
    \frac{\partial (\eta \xi)}{\partial g_i} = 
    \frac{\partial \eta}{\partial g_i} \cdot \epsilon_{\Z_{p}[[F_{n}]]} (\xi) + 
    \eta \cdot \frac{\partial \xi}{\partial g_i} 
    \quad (\eta, \ \xi \in \Z_{p}[[F_{n}]]), 
    $ \\ 
    where $\epsilon_{\Z_{p}[[F_{n}]]} \: \Z_{p}[[F_{n}]] \to \Z_{p}$ 
    is the augmentation map; \\
    \item
    $\displaystyle
    \frac{\partial (w^{-1}) }{\partial g_i} = 
    -w^{-1} \cdot \frac{\partial w}{\partial g_i} \quad (w \in F_{n}). 
    $ 
\end{enumerate}

Now, let $\overline{\Q_{p}}$ be the algebraic closure of the field of $p$-adic numbers 
and assume $H$ is abelian. 
For a positive integer $\ell \in \Z_{\geq 1}$, and for a given homomorphism 
$\varphi \: F_{n} \to \GL(\ell; \overline{\Q_{p}} )$ 
consider a tensor product homomorphism 
$\widetilde{ \varphi } := 
\varphi \otimes (\psi \circ \pi) \: F_{n} \to 
\GL(\ell; \overline{\Q_p} \otimes_{\Z_p} \Z_p[[H]] )$
defined by 
\begin{equation} \label{eq:tensorRep}
\widetilde{ \varphi } (x) := 
(\varphi \otimes (\psi \circ \pi ) ) (x) := 
\varphi(x) \cdot (\psi \circ \pi )(x) = 
\GL(\ell; \overline{\Q_p} \otimes_{\Z_p} \Z_p[[H]] ) 
\end{equation}
for any $x \in F_{n}$. 
Consider the $m \ell \times n \ell$ matrix
\[
Q_{\widetilde{\varphi}} = 
\left( \widetilde{\varphi} 
\left( \frac{ \partial r_{j} }{ \partial g_{i} } \right) 
\right)
\in {\rm M}( m \ell, n \ell; \overline{\Q_p} \otimes_{\Z_p} \Z_p[[H]] ), 
\]
as in \cref{prop*:presMat}. 
When $\mathbb{Z}_{p}[[H]]$ is a Noetherian UFD, 
for a non-negative integer $d$ with $d < n \ell$, 
we define the {\it $d$-th Fitting ideal} 
$E_{d}( A_{\widetilde{\varphi}} ) \subset \overline{\Q_p} \otimes_{\Z_p} \Z_p [[H]]$ of 
$A_{\widetilde{\varphi}}$ as the ideal generated by 
the $(n \ell - d)$-minors of $Q_{\widetilde{\varphi}}$. 
If $n \ell - d > m \ell$, we let $E_{d}( A_{\widetilde{\varphi}} ) := 0$. 
These ideals depend only on $A_{\widetilde{\varphi}}$ and are independent of the choice of a presentation 
(\cite[Section 3.1, Theorem 1]{northcott1976finite}). 
Since the GCD of the generators of $E_{d}( A_{\widetilde{\varphi}} )$ is well-defined 
up to multiplication by a unit of 
$\overline{\Q_p} \otimes_{\Z_p} \Z_p[[H]]$, 
we denote the GCD of generators of $E_{d}( A_{\widetilde{\varphi}} )$ by 
$\Delta_d( A_{\widetilde{\varphi}} ) \in \overline{\Q_p} \otimes_{\Z_p} \Z_p [[H]]$. 
If $\mathrm{Im} \varphi \subset \mathrm{GL}(\ell,\mathbb Z_p)$, then 
$E_d( A_{\widetilde{\varphi}} )$ can be defined as an ideal of $\mathbb Z_p[[H]]$,
and so $\Delta_d( A_{\widetilde{\varphi}} )$ is the GCD in $\mathbb Z_p[[H]]$. 

\section{Burde--de Rham Theorem for finitely presented pro-$p$ groups} \label{sec:mainThm}

Denote $\Gamma := \mathbb{Z}_{p}$, and 
suppose that we can take a surjective homomorphism 
$\alpha \: G \to \Gamma$ such that 
$\alpha(g_{i}) = \gamma$ 
for all $i$, 
where $\gamma$ is a topological generator of $\Gamma$. 

For an integer $k \in \Z_{\geq 2}$, and 
for a given homomorphism 
$\varphi_{k-2} \: F_{n} \to \GL(k-1; \overline{\Q_{p}} )$, 
consider a tensor product homomorphism 
$\widetilde{ \varphi_{k-2} } := 
\varphi_{k-2} \otimes (\alpha \circ \pi) \: F_{n} \to 
\GL(k-1; \overline{\Q_p} \otimes_{\Z_p} \Z_p[[\Gamma]] )$ as 
\eqref{eq:tensorRep}. 
For a given 
$a_{k-1} \in \overline{\Q_{p}}^{\times}$ 
such that $|a_{k-1}-1|_p<1$, 
let $\overline{\beta_{k-1}} \: F_{n} \to \overline{\Q_{p}}^{\oplus (k-1)}$ 
be a crossed homomorphism with respect to 
$\widetilde{\varphi_{k-2}} |_{\gamma = a_{k-1}} \: F_{n} \to \GL(k-1; \overline{\Q_{p}} )$, 
namely satisfying the relation 
$$\overline{\beta_{k-1}}(x y) = 
\overline{\beta_{k-1}}(x) + 
\widetilde{\varphi_{k-2}}(x) |_{\gamma = a_{k-1}} \overline{\beta_{k-1}}(y) 
= \overline{\beta_{k-1}}(x) + a_{k-1}^{ ( \alpha_{\gamma} \circ \pi ) (x)} \cdot 
\varphi_{k-2}(x) \cdot \overline{\beta_{k-1}}(y)$$ 
for any $x$, $y\in F_{n}$, 
where $\alpha_{\gamma} \: G \to \Z_{p}$ maps 
$x \in G$ to the power of $\gamma$ of $\alpha(x)$. 
Note that by the definition of the complete differential module, 
the crossed homomorphism $\overline{\beta_{k-1}} \: F_{n} \to \overline{\Q_{p}}^{\oplus (k-1)}$ 
can be regarded as 
a $\overline{\Q_{p}}$-module homomorphism from $A_{ \widetilde{\varphi_{k-2}} }$ to 
$\overline{\Q_{p}}^{ \oplus (k-1) }$. 
Namely, we have an isomorphism between the set of a crossed homomorphism with respect to 
$\widetilde{\varphi_{k-2}} |_{\gamma = a_{k-1}}$ 
and the set 
${\rm Hom}_{ \overline{\Q_{p}} }(A_{ \widetilde{\varphi_{k-2}} }, \overline{\Q_{p}}^{\oplus (k-1)} )$ 
as vector spaces over $\overline{\Q_{p}}$. 
Using the crossed homomorphism $\overline{\beta_{k-1}} \: F_{n} \to \overline{\Q_{p}}^{\oplus (k-1)}$, 
we define a homomorphism 
$\varphi_{k-1} \: F_{n} \to \GL(k; \overline{\Q_{p}} )$ as 
$$
\varphi_{k-1}(g_{1}) =
\begin{pmatrix} 
\widetilde{\varphi_{k-2}} (g_{1}) |_{\gamma = a_{k-1}} & \overline{\beta_{k-1}}(g_{1}) \\ 
0 & 1 
\end{pmatrix}, \ 
\cdots, \ 
\varphi_{k-1}(g_{n}) = 
\begin{pmatrix} 
\widetilde{\varphi_{k-2}} (g_{n}) |_{\gamma = a_{k-1}} & \overline{\beta_{k-1}}(g_{n}) \\ 
0 & 1 
\end{pmatrix}. 
$$

We focus on 
when $\varphi_{k-1} \: F_{n} \to \GL(k; \overline{\Q_{p}} )$  
factors through 
$G$ as a representation 
$\rho_{k-1} 
\: G \to \GL(k; \overline{\Q_{p}} )$, 
where $\rho_{k-1} \circ \pi = \varphi_{k-1}$. 
We say such a representation 
$\rho_{k-1} 
\: G \to \GL(k; \overline{\Q_{p}} )$ 
an {\it extension of a homomorphism} 
$\varphi_{k-1} \: F_{n} \to \GL(k; \overline{\Q_{p}} )$, 
or the homomorphism $\varphi_{k-1} \: F_{n} \to \GL(k; \overline{\Q_{p}} )$ 
{\it can be extended to} 
a representation $\rho_{k-1} \: G \to \GL(k; \overline{\Q_{p}} )$. 

We have the following main results, 
which is an analogue of the Burde--de Rham theorem including higher-dimensional cases 
(\cite{kitano2013linear, HeusenerPorti2015rep, HeusenerPortiSuarez2001defoRed, SilverWilliams2011theorem}): 

\begin{thm} \label{thm:main}
Let $G$ be a finitely presented pro-$p$ group with the total degrees of all relators being $0$. 
Let $\varphi_{k-2} \: F_{n} \to \GL(k-1; \overline{\Q_{p}} )$ be 
a homomorphism that factors through 
$G$. 
Then the homomorphism 
defined by 
$$\varphi_{k-1} \: F_{n} 
\ni g_{i} 
\mapsto 
\begin{pmatrix} 
\widetilde{ \varphi_{k-2} } (g_{i})|_{\gamma = a_{k-1}} & 
\overline{\beta_{k-1}}(g_{i}) \\ 
\mathbf{0} & 1 
\end{pmatrix} 
\in \GL(k; \overline{\Q_{p}} )$$ 
for all $i$, 
factors through 
$G$ 
as a representation 
$\rho_{k-1} \: G \to \GL(k; \overline{\Q_{p}} )$ 
if and only if 
$Q_{\widetilde{\varphi_{k-2}}} |_{\gamma = a_{k-1} } \cdot \mathbf{b} = 0$, 
where $\mathbf{b} = {}^t (\overline{\beta_{k-1}}(g_{1}), \ldots, \overline{\beta_{k-1}}(g_{n}))
\in \overline{\Q_{p}}^{\oplus n (k-1)}$. 

Moreover, we have 
$
\dim_{\overline{\Q_{p}}}
{\rm Hom}_{ \overline{\Q_{p}} }(A_{ \widetilde{\varphi_{k-2}} }, \overline{\Q_{p}}^{\oplus (k-1)} ) 
\geq k, 
$
namely, we have at least 
$k$ linearly independent extensions 
$\rho_{k-1} \: G \to \GL(k; \overline{\Q_{p}} )$ 
of a homomorphism 
$\varphi_{k-1} \: F_{n} \to \GL(k; \overline{\Q_{p}} )$ 
if 
and only if 
$a_{k-1} \in \overline{\Q_{p}}^{\times}$ is a zero of 
$\Delta_{k-1}( A_{ \widetilde{\varphi_{k-2}} } ) \in 
\overline{\Q_p} \otimes_{\Z_p} \Z_p[[\Gamma]]$. 
\end{thm}

\noindent
{\it Proof. }
The proof is based on \cite[Section 4]{HeusenerPortiSuarez2001defoRed}. 
We use the following \cref{lem:main}: 

\begin{lem*}{\hspace*{-1.2mm}} \label{lem:main}
For $w = w(g_{1}, \dots , g_{n}) \in F_{n}$, 
we have 
$$
\varphi_{k-1}(w) = 
\begin{pmatrix} 
\widetilde{ \varphi_{k-2} } (w)|_{\gamma = a_{k-1}} & 
\sum\limits_{i=1}^{n} 
\displaystyle \left( \widetilde{ \varphi_{k-2} } 
\left. \left( \frac{\partial w}{ \partial g_{i} } \right) \right) \right|_{\gamma = a_{k-1}}
\cdot 
\overline{\beta_{k-1}}(g_{i}) \\ 
\mathbf{0} & 1 
\end{pmatrix} 
\in \GL(k; \overline{\Q_{p}} ). 
$$
\end{lem*}

\noindent
{\it Proof. }
Due to \cite[Section 2]{ihara1986galois}, 
since the classical Fox derivative 
$\frac{\partial}{ \partial g_{i} } \: \Z[F_{n}] \to \Z[F_{n}]$ 
of group algebras is continuous with respect to the topology of $\Z[F_{n}]$, 
which is induced from the topology of $\Z_{p}[[F_{n}]]$ via the canonical injection $\mathbb{Z}[F_{n}] \to \Z_{p}[[F_{n}]]$, 
it is enough to prove the classical case. 

We denote by 
the same $\pi$ for the algebra homomorphism 
$\mathbb{Z}[F_{n}] \to \mathbb{Z}[G]$ of group algebras 
induced by $\pi \: F_{n} \to G$, and 
the same $\alpha$ for the algebra homomorphism 
$\mathbb{Z}[G] \to \mathbb{Z}[\Gamma]$ of group algebras 
induced by $\alpha \: G \to \Gamma$. 
For simplicity, denote 
$\overline{\eta} := 
(\alpha \circ \pi) (\eta) \in \Z[\Gamma]$, and 
$\partial_{i} \eta := \frac{\partial \eta}{ \partial g_{i} } \in \Z[F_{n}]$ for 
$\eta \in \Z[F_{n}]$. 
We prove this by induction on the dimension $k$ and the length $l$ of the word $w$. 

For $k = 2$ and for $l = 1$, since 
$$
\partial_{i} ( g_{h} ) = 
\begin{cases}
1 & (i = h) \\
0 & (i\neq h)
\end{cases}, \ 
\partial_{i} ( g_{h}^{-1} ) = 
\begin{cases}
-g_{h}^{-1} & (i = h) \\
0 & (i\neq h)
\end{cases}, 
$$
we have 
\begin{align*}
\varphi_{1}(g_{h}) &= 
\begin{pmatrix} 
a & \overline{\beta_{1}}(g_{h}) \\ 
0 & 1 
\end{pmatrix} = 
\begin{pmatrix} 
\overline{g_{h}} |_{\gamma = a} & 
\sum\limits_{i=1}^{n} 
( \overline{\partial_{i} g_{h}} ) |_{\gamma = a} \cdot \overline{\beta_{1}}(g_{i}) \\ 
0 & 1 
\end{pmatrix}, \\ 
\varphi_{1}(g_{h}^{-1}) &=
\begin{pmatrix} 
a^{-1} & -a^{-1} \cdot \overline{\beta_{1}}(g_{h}) \\ 
0 & 1 
\end{pmatrix} = 
\begin{pmatrix} 
\overline{g_{h}^{-1}} |_{\gamma = a} & 
\sum\limits_{i=1}^{n} ( \overline{\partial_{i} g_{h}^{-1} } ) |_{\gamma = a} \cdot \overline{\beta_{1}}(g_{i}) \\ 0 & 1 
\end{pmatrix}, 
\end{align*}
and so the claim holds. 

Next, assume that the claim holds for $k$ and $l-1$. 
Namely, for any $w_{l-1} \in F_{n}$ with the length $l-1$, 
assume we have 
$$
\varphi_{k-1}(w_{l-1}) = 
\begin{pmatrix} 
\widetilde{ \varphi_{k-2} } (w_{l-1})|_{\gamma = a_{k-1}} & 
\sum\limits_{i=1}^{n} 
\displaystyle \widetilde{ \varphi_{k-2} }
\left. \left( \partial_{i} w_{l-1} \right) \right|_{\gamma = a_{k-1}}
\cdot \overline{\beta_{k-1}}(g_{i}) \\ 
\mathbf{0} & 1 
\end{pmatrix} 
\in \GL(k; \overline{\Q_{p}} ). 
$$

\noindent
Then for the case of $k$ and $l$, since 
$$
\partial_{i} ( w_{l-1} g_{h} ) = 
\begin{cases}
\partial_{i} w_{l-1} + w_{l-1} & (i = h) \\
\partial_{i} w_{l-1} & (i\neq h)
\end{cases}, \ 
\partial_{i} ( w_{l-1} g_{h}^{-1} ) = 
\begin{cases}
\partial_{i} w_{l-1} - w_{l-1} g_{h}^{-1} & (i = h) \\
\partial_{i} w_{l-1} & (i\neq h)
\end{cases}, 
$$
we have 
$$\leqno{\varphi_{k-1}(w_{l-1} g_{h}) = \varphi_{k-1}(w_{l-1}) \cdot \varphi_{k-1}(g_{h}) }$$
\begin{align*}
&= \begin{pmatrix} 
\widetilde{ \varphi_{k-2} } (w_{l-1} g_{h})|_{\gamma = a_{k-1}} & 
\sum\limits_{i=1}^{n} 
\displaystyle 
\widetilde{ \varphi_{k-2} } ( \partial_{i} w_{l-1} ) |_{\gamma = a_{k-1}} 
\cdot \overline{\beta_{k-1}}(g_{i}) - 
\widetilde{ \varphi_{k-2} }(w_{l-1}) |_{\gamma = a_{k-1}} 
\cdot \overline{\beta_{k-1}}(g_{h}) \\ 
\mathbf{0} & 1 
\end{pmatrix} \\ 
&= \begin{pmatrix} 
\widetilde{ \varphi_{k-2} } (w_{l-1} g_{h})|_{\gamma = a_{k-1}} & 
\sum\limits_{i=1}^{n} 
\displaystyle 
\widetilde{ \varphi_{k-2} }
\left. \left( \partial_{i} (w_{l-1} g_{h}) \right) \right|_{\gamma = a_{k-1}}
\cdot \overline{\beta_{k-1}}(g_{i}) \\ 
\mathbf{0} & 1 
\end{pmatrix} 
\in \GL(k; \overline{\Q_{p}} ), 
\end{align*}
$$\leqno{\varphi_{k-1}(w_{l-1} g_{h}^{-1}) = 
\varphi_{k-1}(w_{l-1}) \cdot \varphi_{k-1}(g_{h})^{-1} }$$
\begin{align*}
&= \begin{pmatrix} 
\widetilde{ \varphi_{k-2} } (w_{l-1} g_{h}^{-1})|_{\gamma = a_{k-1}} & 
\sum\limits_{i=1}^{n} 
\displaystyle 
\widetilde{ \varphi_{k-2} } ( \partial_{i} w_{l-1} ) |_{\gamma = a_{k-1}} 
\cdot \overline{\beta_{k-1}}(g_{i}) - 
\widetilde{ \varphi_{k-2} }(w_{l-1}) |_{\gamma = a_{k-1}} 
\cdot a^{-1} \cdot \overline{\beta_{k-1}}(g_{h}) \\ 
\mathbf{0} & 1 
\end{pmatrix} \\ 
&= \begin{pmatrix} 
\widetilde{ \varphi_{k-2} } (w_{l-1} g_{h}^{-1})|_{\gamma = a_{k-1}} & 
\sum\limits_{i=1}^{n} 
\displaystyle 
\widetilde{ \varphi_{k-2} }
\left. \left( \partial_{i} (w_{l-1} g_{h}^{-1}) \right) \right|_{\gamma = a_{k-1}} 
\cdot \overline{\beta_{k-1}}(g_{i}) \\ 
\mathbf{0} & 1 
\end{pmatrix} 
\in \GL(k; \overline{\Q_{p}} ), 
\end{align*}
and so the claim holds. 

Lastly, for the case of $k + 1$ and $l-1$, the claim holds 
using a calculation similar to the one above. 
\hfill $\Box$

\ 

By the assumption that 
the total degree $\deg(r_{j})$ of a relator $r_{j}$ of $G$ is $0$ for all $j$, we have 
$$
\widetilde{ \varphi_{k-2} } (r_{j})|_{\gamma = a_{k-1}} = 
a_{k-1}^{\deg(r_{j})} \cdot \varphi_{k-2} (r_{j}) = \varphi_{k-2} (r_{j}). 
$$ 
Hence, by \cref{lem:main}, 
the homomorphism $\varphi_{k-1} \: F_{n} \to \GL(k; \overline{\Q_{p}} )$ 
factors through 
$G$ as a representation 
$\rho_{k-1}\: G \to \GL(k; \overline{\Q_{p}} )$ 
if and only if 
$$
\sum\limits_{i=1}^{n} 
\displaystyle \left( \widetilde{ 
\varphi_{k-2} } 
\left. \left( \frac{\partial r_{j}}{ \partial g_{i} } \right) \right) 
\right|_{\gamma = a_{k-1} } 
\cdot \overline{\beta_{k-1}}(g_{i}) = 0 
$$ 
for all $j = 1, \dots, m$. 
This means that the following linear system has a non-zero solution 
$\mathbf{b}_{k-1} \in \overline{\Q_{p}}^{\oplus n(k-1)}$: 
\begin{equation} \label{eq:linSys}
Q_{ \widetilde{\varphi_{k-2}} } |_{\gamma = a_{k-1}} \cdot \mathbf{b}_{k-1} = 0.     
\end{equation}

\noindent
In order to have at least $k$ 
linearly independent extensions 
of $\rho_{k-1}$, 
we need to have at least $k$ solutions 
of \eqref{eq:linSys}, 
which are linearly independent. 
This condition is equivalent to
$$
{\rm rank} 
( Q_{\widetilde{\varphi_{k-2}}} |_{\gamma = a_{k-1}} ) 
\leq n (k - 1) - k
< n (k - 1) - (k - 1), 
$$
which also 
means that $\Delta_{k - 1}( A_{ \widetilde{ \varphi_{k - 2} } } ) = 0$, and
$a_{k-1} \in \overline{\Q_{p}}^{\times}$ is a zero of 
$\Delta_{k-1}( A_{ \widetilde{\varphi_{k-2}} } ) \in 
{\overline{\Q_p} \otimes_{\Z_p} \Z_p[[ \Gamma ]]}$. 

\hfill $\Box$

\begin{rem} \label{rem:analogy}
{\rm
In the spirit of arithmetic topology, we call 
$\Delta_{k-1}(A_{ \widetilde{ \varphi_{k-2} }}) \in 
\overline{\Q_p} \otimes_{\Z_p} \Z_p[[ \Gamma ]]
$ 
the {\it arithmetic twisted Alexander polynomial associated with} $\rho_{k-2}$ 
which may be regarded as a knot-theoretic analogue of 
the twisted Alexander invariant associated with a representation of a knot group 
(\cite{lin2001representations, wada1994twisted}) or 
the $L$-function 
associated with the universal representation of a knot group 
(\cite{kmtt2018certain, kmtt2022adjoint, MTTU2017universal}). 
These knot invariants would be interesting to compare in our future study. 

It may be interesting to note that 
the knot theoretic analogue of \cref{thm:main} 
for higher-dimensional cases 
are proven by Heusener and Porti 
\cite[Theorem 1.4]{HeusenerPorti2015rep}. 
Their result gives a partially positive answer to 
the rigidity questions of Mazur (\cite[p.452 Proposition, p.454 Question 3]{mazur2000variation}). 
Their idea of using the non-vanishing of the cup products is reminiscent of Sharifi's works 
(\cite{sharifi2007iwasawa, sharifi2007massey}), and 
would also be of interest to us in our future study. 
}
\end{rem}

\section{Examples} \label{sec:Eg}

Let us consider some examples under Iwasawa theoretic conditions. 
Let $K_{\infty}/K$ a $\mathbb Z_p$-extension, and 
$\varSigma$ the set of primes of $K$ ramifying in $K_{\infty}/K$. 
Let $S$ be a finite set consisting of places of $K$, and 
let $(K_{\infty})_S$ (resp.\ $K_S$) be the maximal pro-$p$-extension of 
$K_{\infty}$ (resp.\ $K$) unramified outside $S$. 
Put $G_S(K)=\mathrm{Gal}(K_S/K)$ and $\widetilde{G}_S(K)=\mathrm{Gal}((K_{\infty})_S/K)$. 
It is known 
(\cite[Theorem 10.7.12]{neukirch2020cohomology}, \cite{shafarevich1963ext} 
(resp. \cite{blondeau2013cohomological, ElMizusawa2022pro, mizusawa2019pro, salle2008pro}))
that the Galois group 
$G_{S}(K)$ (resp.\ $\widetilde{G}_{S}(K)$) 
is finitely presented as a pro-$p$ group. 
Moreover, if $\varSigma \subset S$, then we have 
$(K_{\infty})_S=K_S$ and $\widetilde{G}_S(K)=G_S(K)$. 
Hence, we focus on 
the Galois $\widetilde{G}_{S}(K)$ having 
the same presentation as \eqref{eq:grpPres}: 
$$
\widetilde{G}_S(K) =
\langle g_{1}, \dots, g_{n} \mid 
r_1 = \cdots = r_m = 1 \rangle^{{\rm pro-}p}. 
$$
The total degrees 
$\deg(r_{j})$ of all relators $r_{j}$ of 
$\widetilde{G}_S(K)$ 
are always $0$ 
since 
$\alpha(r_{j}) = \gamma^{\deg(r_{j})}$ for all $j$, and
the ramified pro-$p$ extension 
$(K_{\infty})_{S} / K$ 
contains the $\Z_{p}$-extension 
$K_{\infty} / K$, 
which means 
the quotient $\Z_{p}$ of $\widetilde{G}_S(K)$ does not vanish. 
Since $S$ contains the set of primes of $K$ over $p$, 
we can take 
the profinite group 
$H$ in \cref{sec:DiffMod} to be 
the Galois group 
$\mathrm{Gal}(K_{\infty} / K) \simeq \mathbb{Z}_{p} = \Gamma$. 
Moreover, the complete $\psi$-differential module 
$A_{\psi}$ becomes a module over 
$\Lambda := \mathbb{Z}_{p}[[ \Gamma ]] \simeq 
\mathbb{Z}_{p}[[ T ]]$ by 
the pro-$p$ Magnus isomorphism 
(\cite[Lemma 9.3.1]{morishita2024knots}). 
Since $\Gamma$ 
is an abelian group, and 
$\Lambda$ is a Noetherian UFD, the GCD $\Delta_{d}(A_{\psi})$
of the $d$-th Fitting ideal 
$(d \geq 0)$ 
of $A_{\psi}$ 
is defined, and we can regard 
$\Delta_{d}(A_{\psi})$ as 
an element in 
$
\overline{\Q_p} \otimes_{\Z_p} \Z_p[[ \Gamma ]]
$. 

Denote $N := \Ker(\psi \: \widetilde{G}_S(K) \to \Gamma)$, so that we have a short exact sequence of pro-$p$ groups 
$$
1 \longrightarrow N \longrightarrow 
\widetilde{G}_S(K) 
\stackrel{\psi}{\longrightarrow} \Gamma \longrightarrow 1. 
$$
Then we have the exact sequence of left $\Lambda$-modules 
\begin{equation} \label{eq:exactSeq}
0 \longrightarrow N^{\rm ab} \stackrel{\theta_{1}}{\longrightarrow} 
A_{\psi} \stackrel{\theta_{2}}{\longrightarrow} 
\Z_{p}[[ \Gamma ]] \stackrel{\epsilon_{\Z_{p}[[ H ]]} }{\longrightarrow} \Z_{p} \longrightarrow 0,     
\end{equation}
where 
$N^{\rm ab} = H_{1}(N, \Z_{p})$ 
is the maximal pro-$p$ quotient of 
the abelianization $N / [N, N]$ of $N$, 
$\theta_{1}$ is the homomorphism induced by $n \mapsto d n$ for any $n \in N$, and 
$\theta_{2}$ is the homomorphism induced by $dg \mapsto \psi(g) - 1$ for any 
$g \in \widetilde{G}_S(K)$ 
(\cite[Section 10.4]{morishita2024knots}). 
We call the left $\Lambda$-module $N^{\rm ab}$ the {\it Iwasawa module}. 
By \eqref{eq:exactSeq}, and 
since the kernel $\Ker( \epsilon_{ \Lambda } \: \Lambda \to \Z_{p} )$ 
of the augmentation map $\epsilon_{ \Lambda } \: \Lambda \to \Z_{p}$ 
is $\Lambda$-isomorphic to $\Lambda$, 
we have a $\Lambda$-isomorphism $A_{\psi} \simeq N^{\rm ab} \oplus \Lambda$. 
Therefore, we have
$\Delta_{d} (N^{\rm ab}) = \Delta_{d+1}(A_{\psi})$
for any $d \geq 0$. 

Here are some examples of extensions of homomorphism 
for finitely presented pro-$p$ Galois groups $G_{S}(K)$ 
over $\Z_{p}$-extensions. 
For simplicity, 
we consider the case where $\varphi_{0}$ is the constant map to $1$, so that we have 
$\widetilde{\varphi_{0}} = \alpha \circ \pi$. 
We used 
\texttt{SageMath} \cite{sagemath} for calculations of group presentations and Fox derivatives, and 
\texttt{Mathematica} \cite{Mathematica} for calculations of linear algebra. 

\begin{eg} \label{eg:GL2}
{\rm
Let $p$ be an odd prime number, and 
$S = \{ q_{1}, q_{2} \}$ a set of prime numbers
$q_{1} \equiv q_{2} \equiv 1 \ ({\rm mod} \ p)$, 
$q_1 \not\equiv 1 \ ({\rm mod} \ p^{2})$, 
$q_2 \not\equiv 1 \ ({\rm mod} \ p^{2})$, and 
$G_{ \emptyset }( ( \Q_{ \{ q_{1} \} } \Q_{ \{ q_{2} \} } )_{\infty})$ is trivial. 
Then due to \cite[Theorem 1]{MizusawaOzaki2013tame}, 
the Galois group 
$G_{S}(\Q_{\infty}) = 
\mathrm{Gal}( (\Q_{\infty})_S / \Q_{\infty})$ 
has a presentation 
$$
G_{S}(\Q_{\infty}) = 
\langle x, \ y \mid
x^{p^2} = 1, \ 
y^{-1} x y = x^{1+p} 
\rangle^{{\rm pro-}p}. 
$$
Then  
the presentation of the Galois group 
$\widetilde{G}_S(\Q) = \mathrm{Gal}(( \Q_{\infty})_S / \Q)$ 
can be described as follows:
$$
\widetilde{G}_S(\Q) = 
\left\langle g_1, \ g_2,  \ g_3 \  \middle| \ 
\begin{aligned} 
(g_2 g_1^{-1})^{p^2} = 1, \ 
(g_3 g_1^{-1})^{-1} g_2 g_1^{-1} g_3 g_1^{-1} = (g_2 g_1^{-1})^{1+p}, \\ 
g_1 g_2 = g_2 g_1, \ 
g_1 (g_3 g_1^{-1}) g_1^{-1} = (g_3 g_1^{-1})^{1+p} (g_2 g_1^{-1})^{p} 
\end{aligned}
\right\rangle^{{\rm pro-}p}, 
$$
where 
$\alpha(g_i) = \gamma$ for all $i$. 
Then the matrix 
$Q_{\widetilde{ \varphi_{0} } }$ 
is given by 
$$
Q_{\widetilde{ \varphi_{0} } } 
= 
\begin{pmatrix} 
-p^{2} & p^{2} & 0 \\ 
p & -p & 0 \\ 
\gamma-1 & -(\gamma-1) & 0 \\ 
- \{ \gamma - (1 + 2p) \} & -p & \gamma - (1 + p) \\ 
\end{pmatrix}, 
$$
and the GCD $\Delta_{1}(A_{ \widetilde{\varphi_{0}} })$ 
of all 2-minors of $Q_{ \alpha }$ is given by 
$$
\Delta_{1}(A_{ \widetilde{\varphi_{0}} }) 
\doteq \gamma - (1+p) 
\in \Z_p[[\Gamma]]. 
$$

\noindent
Let $a_{1} = 1 + p 
\in \overline{\Q_{p}}^{\times}$ 
be a zero of 
$\Delta_{1}(A_{ \widetilde{\varphi_{0}} })$. 
Then 
\begin{align*}
Q_{\widetilde{ \varphi_{0} } } |_{\gamma = 1+p} 
=
\begin{pmatrix} 
-p^{2} & p^{2} & 0 \\ 
p & -p & 0 \\ 
p & -p & 0 \\ 
p & -p & 0 \\ 
\end{pmatrix}, 
\end{align*}
and hence $\mathbf{b} = {^t (b_1, b_2, b_3)}$ is a solution of 
$Q_{\widetilde{ \varphi_{0} } } |_{\gamma = 1+p} \cdot \mathbf{b}=\boldsymbol{0}$ if and only if $b_1 = b_2$. 
Therefore, the homomorphism 
$\varphi_{1} \: F_{3} \to \GL(2; \overline{\Q_{p}} )$ defined by 
$$
\varphi_{1}(g_{1}) =
\varphi_{1}(g_{2}) =
\begin{pmatrix} 
1+p & b_{1} \\ 
0 & 1 
\end{pmatrix}, \ 
\varphi_{1}(g_{3}) = 
\begin{pmatrix} 
1+p & b_{3} \\ 
0 & 1 
\end{pmatrix}
$$
can be extended to a representation 
$\rho_{1} \: \widetilde{G}_{S}(\Q) \to \GL(2; \overline{\Q_{p}})$, 
and $\rho_{1}$ is non-abelian
if and only if 
$b_1 \neq b_3$. 
}
\end{eg}

\begin{eg}
{\rm
Let $K_{\infty} / \mathbb Q(\sqrt{-q_1q_2})$ be the cyclotomic $\mathbb Z_2$-extension of 
$\mathbb Q(\sqrt{-q_1q_2})$ with prime numbers $q_1 \equiv 3 \pmod{8}$ and $q_2 \equiv 7 \pmod{16}$. 
By \cite{mizusawa2010maximal}, 
the Galois group 
$N:=\mathrm{Gal}((K_{\infty})_{\emptyset}/K_{\infty})$ is generated by $x, \ y, \ z$ with three 
relators 
\[
r_1=x^2[x,y], \quad 
r_2=x^{-2}[y,z], \quad 
r_3=[x,z], 
\]
and 
$\widetilde{G}_{\emptyset}(\mathbb Q(\sqrt{-q_1 q_2}) ) = 
\mathrm{Gal}( (\mathbb Q(\sqrt{-q_1 q_2}))_{\infty})_{\emptyset} / \mathbb Q(\sqrt{-q_1 q_2}) ) = 
N \rtimes \Gamma$ is generated by $x,y,z,\gamma$ with six relators 
$r_1,r_2,r_3$,
\[
r_4=\gamma x \gamma^{-1}x^{-1}, \quad 
r_5=x^{C_1}y^{-C_0}z^{1-C_1}\gamma z^{-1} \gamma^{-1}, \quad 
r_6=\gamma y^{-1} \gamma^{-1} y z , 
\]
with some $C_0, \ C_1 \in 2\mathbb Z_2$. 
Put $g_1=\gamma$, $g_2=\gamma x$, $g_3=\gamma y$. 
Since $z=(\gamma y^{-1} \gamma^{-1}y)^{-1}$, the Galois group 
$\widetilde{G}_{\emptyset}(\mathbb Q( \sqrt{-q_1q_2} ))$ is generated by $g_1,g_2,g_3$ with five relators 
$r_1,r_2,r_3,r_4,r_5$. 
Then we have 
\begin{align*}
r_1&=g_1^{-1} g_2 g_3^{-1} g_2 g_1^{-1} g_3, \\
r_2&=(g_2^{-1} g_1)^2 g_3^{-1} g_1^2 g_3^{-1} g_1^{-1} g_3^2 g_1^{-1}, \\
r_3&=g_2^{-1} g_1^2 g_3^{-1} g_1^{-1} g_3 g_1^{-1} g_2 g_3^{-1} g_1 g_3 g_1^{-1}, \\
r_4&=g_2 g_1^{-1} g_2^{-1} g_1, \\
r_5&=(g_1^{-1} g_2)^{C_1}(g_1^{-1} g_3)^{-C_0}(g_3^{-1} g_1 g_3 g_1^{-1})^{1-C_1}g_1^2 g_3^{-1} g_1^{-1} g_3 g_1^{-1}, 
\end{align*}
and 
\begin{align*}
Q_{\widetilde{ \varphi_{0} } } 
=
\begin{pmatrix}
-2\gamma^{-1} & 2\gamma^{-1} & 0 \\
2\gamma^{-1} & -2\gamma^{-1} & 0 \\
0 & 0 & 0 \\
-1 + \gamma^{-1} & 1 - \gamma^{-1} & 0 \\
\gamma + C_1 - 2 + (C_0 - 2C_1 + 1)\gamma^{-1} & C_1 \gamma^{-1} & -\gamma - C_1 + 2 +(- C_0 + C_1 - 1)\gamma^{-1} 
\end{pmatrix} . 
\end{align*}
The non-zero $2$-minors 
of $Q_{\widetilde{ \varphi_{0} } }$ are essentially  
\[
2\gamma^{-1}(-\gamma - C_1 + 2 +(- C_0 + C_1 - 1)\gamma^{-1})
\]
and
\[
(1 - \gamma^{-1})(-\gamma - C_1 + 2 +(- C_0 + C_1 - 1)\gamma^{-1}) . 
\]
Their GCD 
$\Delta_{1}(A_{ \widetilde{\varphi_{0}} })$ 
is given by 
\begin{align*}
\Delta_{1}(A_{ \widetilde{\varphi_{0}} }) 
=\gamma^2 +(C_1 - 2)\gamma +(C_0 - C_1 + 1)
\in \Z_p[[\Gamma]]. 
\end{align*}
Let 
$a_{1} = a \in \overline{\Q_{p}}^{\times}$ 
be a zero of 
$\Delta_{1}(A_{ \widetilde{\varphi_{0}} })$. 
Then 
\begin{align*}
Q_{\widetilde{ \varphi_{0} } } |_{\gamma = a} =
\begin{pmatrix}
-2a^{-1} & 2a^{-1} & 0 \\
2a^{-1} & -2a^{-1} & 0 \\
0 & 0 & 0 \\
-1 + a^{-1} & 1 - a^{-1} & 0 \\
-C_1 a^{-1} & C_1 a^{-1} & 0
\end{pmatrix}, 
\end{align*}
and hence $\mathbf{b}={^t(b_1,b_2,b_3)}$ is a solution of 
$Q_{\widetilde{ \varphi_{0} } } |_{\gamma = a} \cdot \mathbf{b}=\boldsymbol{0}$ if and only if $b_1=b_2$. 
Therefore, the homomorphism 
$\varphi_{1} \: F_{3} \to \GL(2; \overline{\Q_{p}} )$ defined by 
$$
\varphi_{1}(g_{1}) =
\varphi_{1}(g_{2}) =
\begin{pmatrix} 
a & b_{1} \\ 
0 & 1 
\end{pmatrix}, \ 
\varphi_{1}(g_{3}) = 
\begin{pmatrix} 
a & b_{3} \\ 
0 & 1 
\end{pmatrix}, 
$$
can be extended to a representation 
$\rho_{1} \: \widetilde{G}_{\emptyset}(\mathbb Q(\sqrt{-q_1q_2})) \to 
\GL(2; \overline{\Q_{p}} )$, 
and $\rho_{1}$ is non-abelian 
if and only if $a \neq 1$ and $b_1 \neq b_3$. 
}
\end{eg}

\begin{eg}
{\rm
Suppose $p \neq 2$, and let $K$ be a real quadratic field satisfying the assumption of 
\cite[Theorem]{komatsu1989maximal}. 
Then 
the Galois group 
$G_{\{p\}}(K)$ 
is generated by $x, \ y$ with a single relator 
\[
r = (y^{d} x y^{-d})(y^{d-1} x y^{-(d-1)})^{C_{d-1}} \cdots (y x y^{-1})^{C_{1}} x h
\]
with some $h \in [N,N]$, where 
$N := \Ker(\alpha \: G_{\{p\}}(K) \rightarrow 
\Gamma \ ; \  x \mapsto 1, \ y \mapsto \gamma)$, 
and with some $C_{i} \in \Z_p$, where $C_{i} \equiv (-1)^{d-i} 
\begin{pmatrix} 
d \\ 
i 
\end{pmatrix} \ 
({\rm mod} \ p)
$. 
Put $g_1=y$, $g_2=yx$. 
Then $x=g_1^{-1}g_2$, and 
\begin{align*}
(\alpha \circ \pi)\left(\frac{\partial (y^jxy^{-j})}{\partial g_i}\right)
&= (\alpha \circ \pi)\left(\frac{\partial y}{\partial g_i}\right)
+
\gamma (\alpha \circ \pi)\left(\frac{\partial (y^{j-1}xy^{-j})}{\partial g_i}\right) \\
&= (\alpha \circ \pi)\left(\frac{\partial y}{\partial g_i}\right)
+
\gamma \left( 
(\alpha \circ \pi)\left(\frac{\partial (y^{j-1}xy^{-(j-1)})}{\partial g_i}\right)
+
(\alpha \circ \pi)\left(\frac{\partial y^{-1}}{\partial g_i}\right)
\right)\\
&= (\alpha \circ \pi)\left(\frac{\partial y}{\partial g_i}\right)
+
\gamma \left( 
(\alpha \circ \pi)\left(\frac{\partial (y^{j-1}xy^{-(j-1)})}{\partial g_i}\right)
-\gamma^{-1}
(\alpha \circ \pi)\left(\frac{\partial y}{\partial g_i}\right)
\right)\\
&= 
\gamma 
(\alpha \circ \pi)\left(\frac{\partial (y^{j-1}xy^{-(j-1)})}{\partial g_i}\right) \\
&= \left\{\begin{array}{ll} -\gamma^{-1+j} & (i=1) \\ \gamma^{-1+j} & (i=2) \end{array}\right. . 
\end{align*}
Since $(\alpha \circ \pi)\left(\frac{\partial h}{\partial g_i}\right)=0$, 
\begin{align*}
(\alpha \circ \pi)\left(\frac{\partial r}{\partial g_i}\right)
=\sum_{j=0}^d C_{j} \cdot (\alpha \circ \pi)\left(\frac{\partial (y^jxy^{-j})}{\partial g_i}\right) =(-1)^{i+1}\gamma^{-1}\sum_{j=0}^d C_{j} \gamma^j 
\end{align*}
where $C_d=1$. 
The GCD $\Delta_{1}(A_{ \widetilde{\varphi_{0}} })$ of 
$1$-minors of $Q_{\widetilde{ \varphi_{0} } }$ is given by 
\[
\Delta_{1}(A_{ \widetilde{\varphi_{0}} }) = 
\sum_{j=0}^d C_{j} \gamma^j 
\in \Z_p[[\Gamma]]. 
\]
Let $a_{1} = a \in \overline{\Q_{p}}^{\times}$ be a zero of $\Delta_{1}(A_{ \widetilde{\varphi_{0}} })$. 
Then 
$$Q_{\widetilde{ \varphi_{0} } } |_{\gamma = a} = 
\begin{pmatrix} 
0 & 0 
\end{pmatrix}
$$ 
and hence the homomorphism 
$\varphi_{1} \: 
F_{2} \to \GL(2; \overline{\Q_{p}} )$ defined by 
$$
\varphi_{1}(g_{1}) =
\begin{pmatrix} 
a & b_{1} \\ 
0 & 1 
\end{pmatrix}, \ 
\varphi_{1}(g_{2}) = 
\begin{pmatrix} 
a & b_{2} \\ 
0 & 1 
\end{pmatrix}
$$
for any $b_1$ and $b_2$ 
can be extended to 
a representation 
$\rho_{1} \: G_{\{p\}}(K) \to 
\GL(2; \overline{\Q_{p}} )$, 
and $\rho_{1}$ is non-abelian 
if and only if $a \neq 1$ and $b_1 \neq b_2$. 
}
\end{eg}

Next, we consider 
extensions of homomorphisms to $\GL(3; \overline{\Q_{p}})$ 
using \cref{eg:GL2}. 

\begin{eg} \label{eg:AbelGL3}
{\rm
For the same Galois group $\widetilde{G}_{S}(\Q)$ and 
the representation $\rho_{1} \: \widetilde{G}_{S}(\Q) \to \GL(2; \Q_{p} )$ 
in \cref{eg:GL2}, 
assume $b_1 = b_3 = b$. 
Namely, we consider extensions of the homomorphism 
$\varphi_{2} \: F_{3} \to \GL(3; \overline{ \Q_{p}} )$. 
Consider 
the homomorphism  
$\widetilde{\varphi_{1}} \: F_{3} \to \GL(2; \overline{\Q_p} \otimes_{\Z_p} \Z_p[[\Gamma]] )$
defined by
$$
\widetilde{\varphi_{1}} (g_{1}) =
\widetilde{\varphi_{1}} (g_{2}) =
\widetilde{\varphi_{1}} (g_{3}) = 
\begin{pmatrix} 
(1+p) \cdot \gamma & b \cdot \gamma \\ 
0 & \gamma 
\end{pmatrix}. 
$$
Then the matrix $Q_{ \widetilde{ \varphi_{1} } }$ is given by 
$$
Q_{\widetilde{ \varphi_{1} } } = 
\begin{pmatrix} 
 -p^{2} & 0 & p^{2} & 0 & 0 & 0 \\
 0 & -p^{2} & 0 & p^{2} & 0 & 0 \\
 p & 0 & -p & 0 & 0 & 0 \\
 0 & p & 0 & -p & 0 & 0 \\
 (1 + p) \gamma - 1 & b \gamma & - (1 + p) \gamma + 1 & -b \gamma & 0 & 0 \\
 0 & \gamma-1 & 0 & -\gamma + 1 & 0 & 0 \\
 - (1 + p) \gamma + (1 + 2p) & 
 -b \gamma & -p & 0 & (1 + p) (\gamma-1) & 
 b \gamma \\
 0 & - \gamma + (1 + 2p) & 0 & -p & 0 & \gamma - (1+p) \\
\end{pmatrix}, 
$$
and the GCD $\Delta_{2}(A_{\widetilde{ \varphi_{1} } })$ 
of all $(6-2)$-minors of $Q_{\widetilde{ \varphi_{1} } }$ is given by 
$$
\Delta_{2}(A_{\widetilde{ \varphi_{1} } }) \doteq 
(\gamma - 1) \{ \gamma - (1 + p) \} 
\in \Z_p[[\Gamma]]. 
$$

(i) 
Let $a_{2} = 1$. 
Then 
$$
Q_{\widetilde{ \varphi_{1} } } |_{\gamma = 1} = 
\begin{pmatrix} 
 -p^{2} & 0 & p^{2} & 0 & 0 & 0 \\
 0 & -p^{2} & 0 & p^{2} & 0 & 0 \\
 p & 0 & -p & 0 & 0 & 0 \\
 0 & p & 0 & -p & 0 & 0 \\
 p & b & -p & -b & 0 & 0 \\
 0 & 0 & 0 & 0 & 0 & 0 \\
 p & -b & -p & 0 & 0 & 
 b \\
 0 & 2p & 0 & -p & 0 & -p \\
\end{pmatrix}, 
$$
and hence $\mathbf{b} = {^t (b_{11}, b_{12}, b_{21}, b_{22}, b_{31}, b_{32} )}$ is a solution of 
$Q_{\widetilde{ \varphi_{1} } } |_{\gamma = 1} \cdot \mathbf{b} = \boldsymbol{0}$ if and only if  
$b_{11} = b_{21}$, and 
$b_{12} = b_{22} = b_{32}$, 
namely, 
$$
\mathbf{b} = 
b_{11} \cdot 
\begin{pmatrix} 
1 \\ 0 \\ 1 \\ 0 \\ 1 \\ 0 
\end{pmatrix} + 
b_{12} \cdot 
\begin{pmatrix} 
0 \\ 1 \\ 0 \\ 1 \\ 0 \\ 1 
\end{pmatrix} + 
b_{31} \cdot 
\begin{pmatrix} 
0 \\ 0 \\ 0 \\ 0 \\ 1 \\ 0 
\end{pmatrix}. 
$$
Therefore, the homomorphism 
$\varphi_{2} \: F_{3} \to \GL(3; \overline{ \Q_{p}} )$ 
defined by 
$$
\varphi_{2}(g_{1}) = 
\varphi_{2}(g_{2}) =
\begin{pmatrix} 
1+p & b & b_{11} \\ 
0 & 1 & b_{12} \\ 
0 & 0 & 1 
\end{pmatrix}, \ 
\varphi_{2}(g_{3}) = 
\begin{pmatrix} 
1+p & b & b_{31} \\ 
0 & 1 & b_{12} \\ 
0 & 0 & 1 
\end{pmatrix}, 
$$
can be extended to 
a representation 
$\rho_{2} \: \widetilde{G}_{S}(\Q) \to \GL(3; \overline{\Q_{p}} )$, and 
$\rho_{2}$ is non-abelian if and only if $b_{11} \neq b_{31}$. 

(ii)
Let $a_{2} = 1 + p$. 
Then 
$$
Q_{\widetilde{ \varphi_{1} } } |_{\gamma = 1 + p} = 
\begin{pmatrix} 
 -p^{2} & 0 & p^{2} & 0 & 0 & 0 \\
 0 & -p^{2} & 0 & p^{2} & 0 & 0 \\
 p & 0 & -p & 0 & 0 & 0 \\
 0 & p & 0 & -p & 0 & 0 \\
 p (2 + p) & (1 + p) b & - p (2 + p) & - (1 + p) b & 0 & 0 \\
 0 & p & 0 & -p & 0 & 0 \\
 - p^{2} & 
 - (1 + p) b & -p & 0 & p (1 + p) & 
 (1 + p) b \\
 0 & p & 0 & -p & 0 & 0 \\
\end{pmatrix}, 
$$
and hence $\mathbf{b} = {^t (b_{11}, b_{12}, b_{21}, b_{22}, b_{31}, b_{32} )}$ is a solution of 
$Q_{\widetilde{ \varphi_{1} } } |_{\gamma = 1+p} \cdot \mathbf{b} = \boldsymbol{0}$ for 
$b_{11} = b_{21}$, 
$b_{12} = b_{22}$, and 
$b_{31} = 
b_{11} 
+ \frac{1}{p} \cdot b \cdot b_{12} 
- \frac{1}{p} \cdot b \cdot b_{32}
$, 
namely, 
$$
\mathbf{b} = 
b_{11} \cdot 
\begin{pmatrix} 
1 \\ 0 \\ 1 \\ 0 \\ 1 \\ 0 
\end{pmatrix} + 
b_{12} \cdot 
\begin{pmatrix} 
0 \\ 1 \\ 0 \\ 1 \\ \frac{1}{p} \cdot b \\ 0 
\end{pmatrix} + 
b_{32} \cdot 
\begin{pmatrix} 
0 \\ 0 \\ 0 \\ 0 \\ -\frac{1}{p} \cdot b \\ 1 
\end{pmatrix}. 
$$
Therefore, the homomorphism 
$\varphi_{2} \: 
F_{3} \to \GL(3; \overline{\Q_{p}} )$ defined by 
$$
\varphi_{2}(g_{1}) =
\varphi_{2}(g_{2}) =
\begin{pmatrix} 
(1+p)^2 & (1+p) \cdot b & b_{11} \\ 
0 & 1+p & b_{12} \\ 
0 & 0 & 1 
\end{pmatrix}, \ 
\varphi_{2}(g_{3}) = 
\begin{pmatrix} 
(1+p)^2 & (1+p) \cdot b & b_{31} \\ 
0 & 1+p & b_{32} \\ 
0 & 0 & 1 
\end{pmatrix}, 
$$
can be extended to a representation 
$\rho_{2} \: \widetilde{G}_{S}(\Q) \to \GL(3; \overline{\Q_{p}} )$, 
and $\rho_{2}$ is non-abelian
if and only if 
$b_{11} \neq b_{31}$ or $b_{12} \neq b_{32}$. 
}
\end{eg}

\begin{eg} \label{eg:nonAbelGL3}
{\rm
For the same Galois group $\widetilde{G}_{S}(\Q)$ and 
the representation $\rho_{1} \: \widetilde{G}_{S}(\Q) \to \GL(2; \overline{\Q_{p}} )$
in \cref{eg:GL2}, 
consider the homomorphism 
$\widetilde{\varphi_{1}} \: F_{3} \to \GL(2; \overline{\Q_p} \otimes_{\Z_p} \Z_p[[\Gamma]] )$
defined by 
$$
\widetilde{\varphi_{1}} (g_{1}) =
\widetilde{\varphi_{1}} (g_{2}) =
\begin{pmatrix} 
(1+p) \cdot \gamma & b_{1} \cdot \gamma \\ 
0 & \gamma 
\end{pmatrix}, \ 
\widetilde{\varphi_{1}} (g_{3}) = 
\begin{pmatrix} 
(1+p) \cdot \gamma & b_{3} \cdot \gamma \\ 
0 & \gamma 
\end{pmatrix}. 
$$
Then the matrix $Q_{\widetilde{ \varphi_{1} } }$ is given by 
{\scriptsize
$$\hspace{-1cm}
Q_{\widetilde{ \varphi_{1} } } = 
\begin{pmatrix} 
 -p^{2} & 0 & p^{2} & 0 & 0 & 0 \\
 0 & -p^{2} & 0 & p^{2} & 0 & 0 \\
 p & - (b_1 - b_3) & -p & b_1 - b_3 & 0 & 0 \\
 0 & p & 0 & -p & 0 & 0 \\
 (1 + p) \gamma - 1 & b_1 \gamma & - (1 + p) \gamma + 1 & -b_1 \gamma & 0 & 0 \\
 0 & \gamma-1 & 0 & -\gamma + 1 & 0 & 0 \\
 - (1 + p) \gamma + (1 + 2p) & 
 \{ p (b_1 - b_3) - b_3 \} \gamma - \frac{1}{2} p (9 + p) (b_1 - b_3) & -p & 
 4p (b_1 - b_3) & (1 + p) (\gamma-1) & b_1 \gamma + \frac{1}{2} p (1 + p) (b_1 - b_3) \\
 0 & - \gamma + (1 + 2p) & 0 & -p & 0 & \gamma - (1+p) \\
\end{pmatrix}, 
$$
}and the GCD $\Delta_{2}(A_{ \widetilde{ \varphi_{1} }  })$ 
of all $(6-2)$-minors of $Q_{ \widetilde{ \varphi_{1} } }$ is given by 
\begin{align*}
\Delta_{2}(A_{\widetilde{ \varphi_{1} } }) = 
\gamma - (1 + p) \in \Z_p[[\Gamma]]. 
\end{align*}

Let $a_{2} = 1 + p$, 
which is also a zero of 
$\Delta_{1}(A_{\widetilde{ \varphi_{0} }
}) = \gamma - (1 + p) 
\in \Z_p[[\Gamma]]$ 
(cf. Example 1). 
Then 
{\footnotesize
$$\hspace{-1cm}
Q_{\widetilde{ \varphi_{1} } } |_{\gamma = 1 + p} = 
\begin{pmatrix} 
 -p^{2} & 0 & p^{2} & 0 & 0 & 0 \\
 0 & -p^{2} & 0 & p^{2} & 0 & 0 \\
 p & - (b_1 - b_3) & -p & b_1 - b_3 & 0 & 0 \\
 0 & p & 0 & -p & 0 & 0 \\
 p (2 + p) & (1 + p) b_1 & - p (2 + p) & - (1 + p) b_1 & 0 & 0 \\
 0 & p & 0 & -p & 0 & 0 \\
 -p^2 & - (1 + p) b_3 + \frac{1}{2} p (-7 + p) (b_1 - b_3) &  -p & 
 4p (b_1 - b_3) & p (1 + p) & \frac{1}{2} (1 + p) \{ 2 b_1 + p (b_1 - b_3) \} \\
 0 & p & 0 & -p & 0 & 0 \\
\end{pmatrix}, 
$$
}and hence $\mathbf{b} = {^t (b_{11}, b_{12}, b_{21}, b_{22}, b_{31}, b_{32} )}$ is a solution of 
$Q_{\widetilde{ \varphi_{1} } } |_{\gamma = 1+p } \cdot \mathbf{b} = \boldsymbol{0}$ if and only if 
$b_{11} = b_{21}$, 
$b_{12} = b_{22}$, and 
$b_{31} = 
b_{11} + 
\{ \frac{1}{p} \cdot b_{3} - \frac{1}{2} (b_{1} - b_{3}) \} \cdot b_{12} 
+ \{ - \frac{1}{p} \cdot b_{1} - \frac{1}{2} (b_{1} - b_{3}) \} \cdot b_{32}
$, 
namely, 
$$
\mathbf{b} = 
b_{11} \cdot 
\begin{pmatrix} 
1 \\ 0 \\ 1 \\ 0 \\ 1 \\ 0 
\end{pmatrix} + 
b_{12} \cdot 
\begin{pmatrix} 
0 \\ 1 \\ 0 \\ 1 \\ \frac{1}{p} \cdot b_{3} - \frac{1}{2} (b_{1} - b_{3}) \\ 0 
\end{pmatrix} + 
b_{32} \cdot 
\begin{pmatrix} 
0 \\ 0 \\ 0 \\ 0 \\ -\frac{1}{p} \cdot b_{1} - \frac{1}{2} (b_{1} - b_{3}) \\ 1 
\end{pmatrix}. 
$$
Therefore, the homomorphism 
$\varphi_{2} \: 
F_{3} \to \GL(3; \overline{\Q_{p}} )$ defined by 
$$
\varphi_{2}(g_{1}) =
\varphi_{2}(g_{2}) =
\begin{pmatrix} 
(1+p)^2 & (1+p) \cdot b_{1} & b_{11} \\ 
0 & 1+p & b_{12} \\ 
0 & 0 & 1 
\end{pmatrix}, \ 
\varphi_{2}(g_{3}) = 
\begin{pmatrix} 
(1+p)^2 & (1+p) \cdot b_{3} & b_{31} \\ 
0 & 1+p & b_{32} \\ 
0 & 0 & 1 
\end{pmatrix}, 
$$
can be extended to a representation 
$\rho_{2} \: \widetilde{G}_{S}(\Q) \to \GL(3; \overline{\Q_{p}} )$, 
and $\rho_{2}$ is non-abelian if and only if 
$b_{1} \neq b_{3}$, $b_{11} \neq b_{31}$, or $b_{12} \neq b_{32}$. 
}
\end{eg}

\section{Cohomological interpretations} \label{sec:cohInterpret}

In this Section, we focus on 
cohomological interpretations 
of our observations so far. 
Generalizations to 
higher-dimensional cases are due to \cite[Section 4]{HeusenerPorti2015rep}. 

For a given $a_{k-1} \in \overline{\Q_{p}}^{\times}$ 
such that $|a_{k-1}-1|_p<1$, 
note that we have the relation 
$$\overline{\beta_{k-1}}(x y) =
\overline{\beta_{k-1}}(x) + 
\widetilde{\varphi_{k-2}}(x) |_{\gamma = a_{k-1}} \overline{\beta_{k-1}}(y) 
= \overline{\beta_{k-1}}(x) + 
a_{k-1}^{( \alpha_{\gamma} \circ \pi )(x)} \cdot 
\varphi_{k-2}(x) \cdot \overline{\beta_{k-1}}(y)$$ 
for any $x$, $y \in F_{n}$, 
where 
$\overline{\beta_{k-1}} \: F_{n} \to \overline{\Q_{p}}^{\oplus (k-1)}$ 
is the crossed homomorphism, 
$\widetilde{ \varphi_{k-2} } := 
\varphi_{k-2} \otimes (\alpha \circ \pi) \: F_{n} \to 
\GL(k-1; \overline{\Q_p} \otimes_{\Z_p} \Z_p[[\Gamma]] )$ 
is the tensor product homomorphism, and 
$\alpha_{\gamma} \: G \to \Z_{p}$ maps 
$g \in G$ to the power of $\gamma$ of $\alpha(g)$. 
In particular, when 
$\rho_{k-2} \: G \to \GL(k-1; \overline{\Q_{p}})$ is a representation 
such that 
$\rho_{k-2} \circ \pi = \widetilde{\varphi_{k-2}} |_{\gamma = a_{k-1}}$, 
$$
\beta_{k-1}(x y) = 
\beta_{k-1}(x) + 
\rho_{k-2}(x) \beta_{k-1}(y), 
$$ 
for any $x$, $y \in G$, 
where 
$\beta_{k-1} \: G \to \overline{\Q_{p}}^{\oplus (k-1)}$ 
is a crossed homomorphism such that $\beta_{k-1} \circ \pi = \overline{\beta_{k-1}}$. 
We call $\beta_{k-1} \: G \to \overline{\Q_{p}}^{\oplus (k-1)}$ 
a {\it 1-cocycle with coefficients} 
$( \overline{\Q_{p}}^{\oplus (k-1)} )_{\rho_{k-2}}$. 

On the other hand, 
if there exists $v \in \overline{\Q_{p}}^{\oplus (k-1)}$ such that
$$
\overline{\beta_{k-1}}(x) = 
\widetilde{\varphi_{k-2}}(x) |_{\gamma = a_{k-1}} v - v = 
(\widetilde{\varphi_{k-2}}(x) |_{\gamma = a_{k-1}} - I_{k-1}) v
= (a_{k-1}^{(\alpha_{\gamma} \circ \pi)(x)} 
\cdot \varphi_{k-2}(x) - I_{k-1}) \cdot v
$$
for any $x \in F_{n}$, 
where 
$I_{k-1} \in \GL(k-1; \overline{\Q_{p}} )$ 
is the $(k-1) \times (k-1)$ identity matrix, 
namely, 
\begin{equation} \label{eq:cobdry} 
\beta_{k-1}(x) = 
\rho_{k-2}(x) v - v = 
(\rho_{k-2}(x) - I_{k-1}) v,     
\end{equation}
for any $x \in G$, 
we call 
$\beta_{k-1} \: G \to \overline{\Q_{p}}^{\oplus (k-1)}$ 
a {\it 1-coboundary with coefficients} $( \overline{\Q_{p}}^{\oplus (k-1)} )_{\rho_{k-2}}$. 

We respectively denote the set of continuous 1-cocyles and 
continuous 1-coboundaries 
by $Z^{1}(G, ( \overline{\Q_{p}}^{\oplus (k-1)} )_{\rho_{k-2}} )$ and $B^{1}(G, ( \overline{\Q_{p}}^{\oplus (k-1)} )_{\rho_{k-2}} )$. 
We can easily see that 
$B^{1}(G, ( \overline{\Q_{p}}^{\oplus (k-1)} )_{\rho_{k-2}} ) \subset Z^{1}(G, ( \overline{\Q_{p}}^{\oplus (k-1)} )_{\rho_{k-2}} )$, 
and 
we define 
the 
{\it 
first continuous cohomology group 
$H^{1}(G, ( \overline{\Q_{p}}^{\oplus (k-1)} )_{\rho_{k-2}})$ 
with coefficients} $( \overline{\Q_{p}}^{\oplus (k-1)} )_{\rho_{k-2}}$ by 
$$H^{1}(G, ( \overline{\Q_{p}}^{\oplus (k-1)} )_{\rho_{k-2}}) := 
Z^{1}(G, ( \overline{\Q_{p}}^{\oplus (k-1)} )_{\rho_{k-2}} ) / B^{1}(G, ( \overline{\Q_{p}}^{\oplus (k-1)} )_{\rho_{k-2}} ). $$ 
Denote 
the 
{\it first cohomology class of} $\beta_{k-1}$ 
in $H^{1}(G, ( \overline{\Q_{p}}^{\oplus (k-1)} )_{\rho_{k-2}})$ 
by $[ \beta_{k-1} ]$. 

From now on, 
let us consider the relation between the zeros of the arithmetic twisted Alexander polynomial 
$\Delta_{k-1}( A_{ \widetilde{\varphi_{k-2}} } ) \in \overline{\Q_p} \otimes_{\Z_p} \Z_p[[\Gamma]]
$ and 
the non-vanishment of $H^{1}(G, ( \overline{\Q_{p}}^{\oplus (k-1)} )_{\rho_{k-2}})$. 
This may be seen as a generalization of \cite[Corollary 14.4.5]{morishita2024knots}. 

\begin{thm} \label{thm:zeroesCoh}
If $a_{k-1} \in \overline{\Q_{p}}^{\times}$ is a zero of 
$\Delta_{k-1}( A_{ \widetilde{\varphi_{k-2}} } ) \in \overline{\Q_p} \otimes_{\Z_p} \Z_p[[\Gamma]]$, then 
$H^{1}(G, ( \overline{\Q_{p}}^{\oplus (k-1)} )_{\rho_{k-2}})$ does not vanish. 
\end{thm}

\noindent
{\it Proof. }
Since we have an isomorphism 
$Z^{1}(G, ( \overline{\Q_{p}}^{\oplus (k-1)} )_{\rho_{k-2}} ) \simeq 
{\rm Hom}_{ \overline{\Q_{p}} }(A_{ \widetilde{\varphi_{k-2}} }, \overline{\Q_{p}}^{\oplus (k-1)} )$ 
as vector spaces over $\overline{\Q_{p}}$, and 
$
\dim_{\overline{\Q_{p}}}
{\rm Hom}_{ \overline{\Q_{p}} }(A_{ \widetilde{\varphi_{k-2}} }, \overline{\Q_{p}}^{ \oplus (k-1)} ) 
\geq k 
$
by \cref{thm:main}, we have 
$$
\dim_{\overline{\Q_{p}}} 
Z^{1}(G, ( \overline{\Q_{p}}^{\oplus (k-1)} )_{\rho_{k-2}} ) \geq k. 
$$

On the other hand, 
since we have a surjective linear map 
\begin{equation} \label{eq:linMap} 
F \: \overline{\Q_{p}}^{\oplus (k-1)} \to 
B^{1}(G, ( \overline{\Q_{p}}^{\oplus (k-1)} )_{\rho_{k-2}} ), 
\end{equation}
defined by 
$F(v) = \beta_{k-1, v}$, where 
$\beta_{k-1, v} (x) = \rho_{k-2}(x) v - v$ for any $x \in G$, 
we have 
$$
\dim_{\overline{\Q_{p}}} 
B^{1}(G, ( \overline{\Q_{p}}^{\oplus (k-1)} )_{\rho_{k-2}} ) 
\leq k-1. 
$$
Hence, we have 
$H^{1}(G, ( \overline{\Q_{p}}^{\oplus (k-1)} )_{\rho_{k-2}}) \neq 0$. 
\hfill $\Box$ 

\ 

For the converse, we have the following: 

\begin{thm} \label{thm:cohZeroes}
Assume the matrices 
$\{ \rho_{k-2}(x) \}_{x \in G} \subset {\rm GL}(k-1, \overline{\Q_{p}} )$ do not have 
a simultaneous eigenvector with simultaneous eigenvalue $1$. 
If $H^{1}(G, ( \overline{\Q_{p}}^{\oplus (k-1)} )_{\rho_{k-2}})$ does not vanish, then 
$a_{k-1} \in \overline{\Q_{p}}^{\times}$ is a zero of 
$\Delta_{k-1}( A_{ \widetilde{\varphi_{k-2}} } ) \in 
\overline{\Q_p} \otimes_{\Z_p} \Z_p[[\Gamma]]
$. 
\end{thm}

\noindent
{\it Proof. }
Since $H^{1}(G, ( \overline{\Q_{p}}^{\oplus (k-1)} )_{\rho_{k-2}}) \neq 0$, 
we have 
$$
\dim_{\overline{\Q_{p}}} 
B^{1}(G, ( \overline{\Q_{p}}^{\oplus (k-1)} )_{\rho_{k-2}} ) < 
\dim_{\overline{\Q_{p}}} 
Z^{1}(G, ( \overline{\Q_{p}}^{\oplus (k-1)} )_{\rho_{k-2}} ). 
$$
Due to the assumption, 
since we do not have a vector $v \in \overline{\Q_{p}}^{\oplus (k-1)}$ 
such that $\rho_{k-2}(x) v = v$ for any $x \in G$, we have 
$\Ker F = 0$, where 
$F \: \overline{\Q_{p}}^{\oplus (k-1)} \to 
B^{1}(G, ( \overline{\Q_{p}}^{\oplus (k-1)} )_{\rho_{k-2}} )$
is the surjective linear map \eqref{eq:linMap}. 
Hence, we have 
$$
\dim_{\overline{\Q_{p}}} 
B^{1}(G, ( \overline{\Q_{p}}^{\oplus (k-1)} )_{\rho_{k-2}} ) 
= k-1. 
$$
This means that 
$$
\dim_{\overline{\Q_{p}}}
{\rm Hom}_{ \overline{\Q_{p}} }(A_{ \widetilde{\varphi_{k-2}} }, 
\overline{ \Q_{p}}^{\oplus (k-1)} ) = 
\dim_{\overline{\Q_{p}}} 
Z^{1}(G, ( \overline{\Q_{p}}^{\oplus (k-1)} )_{\rho_{k-2}} ) 
\geq k 
$$
and by \cref{thm:main}, 
$a_{k-1} \in \overline{\Q_{p}}^{\times}$ is a zero of 
$\Delta_{k-1}( A_{ \widetilde{\varphi_{k-2}} } ) \in 
\overline{\Q_p} \otimes_{\Z_p} \Z_p[[\Gamma]]$. 
\hfill $\Box$

\

Here are some examples. 
Note that we have an example 
which does not satisfy the assumption of \cref{thm:cohZeroes}
but the claim holds (\cref{eg:cohAbel1GL3}), and
the claim does not hold (\cref{eg:cohAbel1/1+pGL3}). 

\begin{eg} \label{eg:cohGL2}
{\rm
Consider the representation 
$\rho_{1} \: \widetilde{G}_{S}(\Q) \to \GL(2; \overline{\Q_{p}} )$ 
constructed in \cref{eg:GL2}, 
which is defined by 
$$
\rho_{1}(g_{1}) =
\rho_{1}(g_{2}) =
\begin{pmatrix} 
1+p & b_{1} \\ 
0 & 1 
\end{pmatrix}, \ 
\rho_{1}(g_{3}) = 
\begin{pmatrix} 
1+p & b_{3} \\ 
0 & 1 
\end{pmatrix}. 
$$
Let $\beta_{1, (b_{1}, b_{3})} \: 
\widetilde{G}_{S}(\Q)
\to \overline{\Q_{p}}$ 
be a 1-cocycle 
such that 
$\beta_{1, (b_{1}, b_{3})}(g_{1}) = \beta_{1, (b_{1}, b_{3})}(g_{2}) = b_{1}$ and 
$\beta_{1, (b_{1}, b_{3})}(g_{3}) = b_{3}$. 

If $b_{1} = b_{3} = b$, 
namely if $\rho$ is abelian, 
then $\beta_{1, (b, b)}$ is a 1-coboundary 
since we can take $v \in \overline{\Q_{p}}$ as 
$\frac{b}{ (1+p) - 1} = \frac{1}{p} \cdot b$
which satisfies \eqref{eq:cobdry}. 

If $b_{1} \neq b_{3}$, 
namely if $\rho$ is non-abelian, 
then $\beta_{1, (b_{1}, b_{3})}$ is not a 1-coboundary 
since we cannot take $v \in \overline{\Q_{p}}$ which satisfies \eqref{eq:cobdry}. 

Since the cocycles are unique up to adding a coboundary and up to multiplying by a non-zero scalar, 
we have two different cocycles 
$\beta_{1, (1, 1)}$ 
and $\beta_{1, (1, 0)}$, 
where the cohomology class 
$[\beta_{1, (1, 1)}]$ 
is trivial, whereas 
$[\beta_{1, (1, 0)}]$ 
is non-trivial. 
Therefore, we have 
\begin{align*}
\dim_{\overline{\Q_{p}}} Z^{1}(\widetilde{G}_{S}(\Q), ( \overline{\Q_{p}} )_{\rho_{0}} ) &= 2, \\
\dim_{\overline{\Q_{p}}} B^{1}(\widetilde{G}_{S}(\Q), ( \overline{\Q_{p}} )_{\rho_{0}} ) &= 1, \\ 
\dim_{\overline{\Q_{p}}} H^{1}(\widetilde{G}_{S}(\Q), ( \overline{\Q_{p}} )_{\rho_{0}} ) &= 1. 
\end{align*}
}
\end{eg}

\begin{eg} \label{eg:cohAbel1GL3}
{\rm
Consider the representation 
$\rho_{2} \: \widetilde{G}_{S}(\Q) \to \GL(3; \overline{\Q_{p}} )$ 
constructed in \cref{eg:AbelGL3} 
(i), 
which is defined by 
$$
\rho_{2}(g_{1}) =
\rho_{2}(g_{2}) =
\begin{pmatrix} 
1+p & b & b_{11} \\ 
0 & 1 & b_{12} \\ 
0 & 0 & 1 
\end{pmatrix}, \ 
\rho_{2}(g_{3}) = 
\begin{pmatrix} 
1+p & b & b_{31} \\ 
0 & 1 & b_{12} \\ 
0 & 0 & 1 
\end{pmatrix}. 
$$
Let $\beta_{2, (b_{11}, b_{12}, b_{31}, b_{32})} \: \widetilde{G}_{S}(\Q) \to \overline{\Q_{p}}^{\oplus 2}$ 
be a 1-cocycle with coefficient $\varphi_{2}$
such that 
$\beta_{2, (b_{11}, b_{12}, b_{31}, b_{32})} (g_{1}) = 
\beta_{2, (b_{11}, b_{12}, b_{31}, b_{32})} (g_{2}) = {}^{t}(b_{11}, b_{12})$ and 
$\beta_{2, (b_{11}, b_{12}, b_{31}, b_{32})} (g_{3}) = {}^{t}(b_{31}, b_{32})$. 

If $b_{11} = b_{31}$ and $b_{12} = b_{32} = 0$, 
then $\beta_{2, (b_{11}, 0, b_{11}, 0)}$ is a 1-coboundary 
since we can take $v \in \overline{\Q_{p}}^{\oplus 2}$ as 
$
\displaystyle 
\frac{1}{p} 
\begin{pmatrix} 
b_{11} \\ 
0  
\end{pmatrix}
$ which satisfies \eqref{eq:cobdry}. 

For the 
rest of the cases, 
$\beta_{2, (b_{11}, b_{12}, b_{31}, b_{32})}$ is not a 1-coboundary 
since we cannot take $v \in \overline{\Q_{p}}^{\oplus 2}$ which satisfies \eqref{eq:cobdry}. 

Hence, 
we have 
three different cocycles 
$\beta_{2, (1, 0, 1, 0)}$, 
$\beta_{2, (0, 1, 0, 1)}$, and 
$\beta_{2, (0, 0, 1, 0)}$, 
where the cohomology class 
$[\beta_{2, (1, 0, 1, 0)}]$ is trivial, whereas 
$[\beta_{2, (0, 1, 0, 1)}]$ 
and 
$[\beta_{2, (0, 0, 1, 0)}]$
are non-trivial. 
Therefore, we have 
\begin{align*}
\dim_{\overline{\Q_{p}}} Z^{1}(\widetilde{G}_{S}(\Q),  ( \overline{\Q_{p}}^{\oplus 2} )_{\rho_{1}} ) &= 3, \\
\dim_{\overline{\Q_{p}}} B^{1}(\widetilde{G}_{S}(\Q),  ( \overline{\Q_{p}}^{\oplus 2} )_{\rho_{1}} ) &= 1, \\ 
\dim_{\overline{\Q_{p}}} H^{1}(\widetilde{G}_{S}(\Q),  ( \overline{\Q_{p}}^{\oplus 2} )_{\rho_{1}} ) &= 2. 
\end{align*}
}
\end{eg}

\begin{eg} \label{eg:cohAbel1+pGL3}
{\rm
Consider the representation 
$\rho_{2} \: \widetilde{G}_{S}(\Q) \to \GL(3; \overline{\Q_{p}} )$ 
constructed in \cref{eg:AbelGL3} 
(ii), 
which is defined by 
$$
\rho_{2}(g_{1}) =
\rho_{2}(g_{2}) =
\begin{pmatrix} 
(1+p)^2 & (1+p) \cdot b & b_{11} \\ 
0 & 1+p & b_{12} \\ 
0 & 0 & 1 
\end{pmatrix}, \ 
\rho_{2}(g_{3}) = 
\begin{pmatrix} 
(1+p)^2 & (1+p) \cdot b & b_{31} \\ 
0 & 1+p & b_{32} \\ 
0 & 0 & 1 
\end{pmatrix}. 
$$
Let $\beta_{2, (b_{11}, b_{12}, b_{31}, b_{32})} \: \widetilde{G}_{S}(\Q) \to \overline{\Q_{p}}^{\oplus 2}$ 
be a 1-cocycle with coefficient $\varphi_{2}$
such that 
$\beta_{2, (b_{11}, b_{12}, b_{31}, b_{32})} (g_{1}) = 
\beta_{2, (b_{11}, b_{12}, b_{31}, b_{32})} (g_{2}) = {}^{t}(b_{11}, b_{12})$ and 
$\beta_{2, (b_{11}, b_{12}, b_{31}, b_{32})} (g_{3}) = {}^{t}(b_{31}, b_{32})$. 

If $b_{11} = b_{31}$ and $b_{12} = b_{32} = 0$, 
then $\beta_{2, (b_{11}, 0, b_{11}, 0)}$ is a 1-coboundary 
since we can take $v \in \overline{\Q_{p}}^{\oplus 2}$ as 
$
\displaystyle 
\frac{1}{p (2+p)} 
\begin{pmatrix} 
b_{11} \\ 
0  
\end{pmatrix}$ which satisfies \eqref{eq:cobdry}. 

If $b_{11} = b_{31} = 0$ and $b_{12} = b_{32}$, 
then $\beta_{2, (0, b_{12}, 0, b_{12})}$ is a 1-coboundary 
since we can take $v \in \overline{\Q_{p}}^{\oplus 2}$ as 
$
\displaystyle 
\frac{1}{p^{2} (2 + p)} 
\begin{pmatrix} 
- (1 + p) \cdot b \cdot b_{12} \\ 
p (2 + p) \cdot b_{12} 
\end{pmatrix}
$ 
which satisfies \eqref{eq:cobdry}. 

For the rest of the cases, 
$\beta_{2, (b_{11}, b_{12}, b_{31}, b_{32})}$ is not a 1-coboundary 
since we cannot take $v \in \overline{\Q_{p}}^{\oplus 2}$ which satisfies \eqref{eq:cobdry}. 

Hence, we have three different cocycles 
$\beta_{2, (1, 0, 1, 0)}$, 
$\beta_{2, (0, 1, 0, 1)}$, 
and 
$\beta_{2, (0, 0, -\frac{1}{p} \cdot b, 1)}$, 
where the cohomology classes $[\beta_{2, (1, 0, 1, 0)}]$ and 
$[\beta_{2, (0, 1, 0, 1)}]$ 
are trivial, whereas 
$[\beta_{2, (0, 0, -\frac{1}{p} \cdot b, 1)} ]$
is non-trivial. 

Therefore, we have 
\begin{align*}
\dim_{\overline{\Q_{p}}} Z^{1}(\widetilde{G}_{S}(\Q), ( \overline{\Q_{p}}^{\oplus 2} )_{\rho_{1}} ) ) &= 3, \\
\dim_{\overline{\Q_{p}}} B^{1}(\widetilde{G}_{S}(\Q), ( \overline{\Q_{p}}^{\oplus 2} )_{\rho_{1}} ) ) &= 2, \\
\dim_{\overline{\Q_{p}}} H^{1}(\widetilde{G}_{S}(\Q), ( \overline{\Q_{p}}^{\oplus 2} )_{\rho_{1}} ) ) &= 1. 
\end{align*}
}
\end{eg}

\begin{eg} \label{eg:cohNonAbelGL3}
{\rm
Consider the representation 
$\rho_{2} \: \widetilde{G}_{S}(\Q) \to \GL(3; \overline{\Q_{p}} )$ 
constructed in \cref{eg:nonAbelGL3}, 
which is defined by 
$$
\rho_{2}(g_{1}) =
\rho_{2}(g_{2}) =
\begin{pmatrix} 
(1+p)^2 & (1+p) \cdot b_{1} & b_{11} \\ 
0 & 1+p & b_{12} \\ 
0 & 0 & 1 
\end{pmatrix}, \ 
\rho_{2}(g_{3}) = 
\begin{pmatrix} 
(1+p)^2 & (1+p) \cdot b_{3} & b_{31} \\ 
0 & 1+p & b_{32} \\ 
0 & 0 & 1 
\end{pmatrix}. 
$$
Let $\beta_{2, (b_{11}, b_{12}, b_{31}, b_{32})} \: \widetilde{G}_{S}(\Q) \to \overline{\Q_{p}}^{\oplus 2}$ 
be a 1-cocycle with coefficient $\varphi_{2}$
such that 
$\beta_{2, (b_{11}, b_{12}, b_{31}, b_{32})} (g_{1}) = 
\beta_{2, (b_{11}, b_{12}, b_{31}, b_{32})} (g_{2}) = {}^{t}(b_{11}, b_{12})$ and 
$\beta_{2, (b_{11}, b_{12}, b_{31}, b_{32})} (g_{3}) = {}^{t}(b_{31}, b_{32})$. 

If $b_{11} = b_{31}$ and $b_{12} = b_{32} = 0$, 
then $\beta_{2, (b_{11}, 0, b_{11}, 0)}$ is a 1-coboundary 
since we can take $v \in \overline{\Q_{p}}^{\oplus 2}$ as 
$
\displaystyle 
\frac{1}{p (2+p)} 
\begin{pmatrix} 
b_{11} \\ 
0  
\end{pmatrix}$
which satisfies \eqref{eq:cobdry}. 

If $b_{11} = 0$, $b_{12} = b_{32}$, and $b_{31} = -\frac{1+p}{p} \cdot (b_{1} - b_{3}) \cdot b_{12}$, 
then $\beta_{2, (0, b_{12}, -\frac{1+p}{p} \cdot (b_{1} - b_{3}) \cdot b_{12}, b_{12})}$ is a 1-coboundary 
since we can take $v \in \overline{\Q_{p}}^{\oplus 2}$ as 
$
\displaystyle 
\frac{1}{p^{2} (2 + p)} 
\begin{pmatrix} 
- (1 + p) \cdot b_{1} \cdot b_{12} \\ 
p (2 + p) \cdot b_{12} 
\end{pmatrix}
$
which satisfies \eqref{eq:cobdry}. 

For the rest of the cases, 
$\beta_{2, (b_{11}, b_{12}, b_{31}, b_{32})}$ is not a 1-coboundary 
since we cannot take $v \in \overline{\Q_{p}}^{\oplus 2}$ which satisfies \eqref{eq:cobdry}. 

Hence, we have three different cocycles 
$\beta_{2, (1, 0, 1, 0)}$, 
$\beta_{2, (0, 1, -\frac{1+p}{p} \cdot (b_{1} - b_{3}), 1)}$, 
and 
$\beta_{2, (0, 0, -\frac{1}{p} \cdot b_{3} - \frac{1}{2} (b_{1} - b_{3}), 1)}$, 
where the cohomology classes 
$[\beta_{2, (1, 0, 1, 0)}]$ and 
$[\beta_{2, (0, 1, -\frac{1+p}{p} \cdot (b_{1} - b_{3}), 1)}]$ 
are trivial, whereas 
$[\beta_{2, (0, 0, -\frac{1}{p} \cdot b_{1} - \frac{1}{2} (b_{1} - b_{3}), 1)} ]$
is non-trivial. 
Therefore, we have 
\begin{align*}
\dim_{\overline{\Q_{p}}} Z^{1}(\widetilde{G}_{S}(\Q), ( \overline{\Q_{p}}^{\oplus 2} )_{\rho_{1}} ) ) &= 3, \\
\dim_{\overline{\Q_{p}}} B^{1}(\widetilde{G}_{S}(\Q), ( \overline{\Q_{p}}^{\oplus 2} )_{\rho_{1}} ) ) &= 2, \\
\dim_{\overline{\Q_{p}}} H^{1}(\widetilde{G}_{S}(\Q), ( \overline{\Q_{p}}^{\oplus 2} )_{\rho_{1}} ) ) &= 1. 
\end{align*}
}
\end{eg}

\begin{rem}
{\rm
Since 
\begin{align*}
\frac{1}{4} \cdot 
\begin{pmatrix} 
b^2 \\ 2 b \\ b^2 \\ 2 b 
\end{pmatrix} 
&= 
\frac{b^2}{4} \cdot 
\begin{pmatrix} 
1 \\ 0 \\ 1 \\ 0 
\end{pmatrix} + 
\frac{b}{2} \cdot 
\begin{pmatrix} 
0 \\ 1 \\ 0 \\ 1 
\end{pmatrix} + 
0 \cdot 
\begin{pmatrix} 
0 \\ 0 \\ 1 \\ 0 
\end{pmatrix} \\ 
&= 
\frac{b^2}{4} \cdot 
\begin{pmatrix} 
1 \\ 0 \\ 1 \\ 0 
\end{pmatrix} + 
\frac{b}{2} \cdot 
\begin{pmatrix} 
0 \\ 1 \\ 0 \\ 1 
\end{pmatrix} + 
0 \cdot 
\begin{pmatrix} 
0 \\ 0 \\ -\frac{1}{p} \cdot b \\ 1 
\end{pmatrix}, \\ 
\frac{1}{4} \cdot 
\begin{pmatrix} 
b_{1}^2 \\ 2 b_{1} \\ b_{3}^2 \\ 2 b_{3}
\end{pmatrix} 
&= 
\frac{b_{1}^2}{4} \cdot 
\begin{pmatrix} 
1 \\ 0 \\ 1 \\ 0 
\end{pmatrix} + 
\frac{b_{1}}{2} \cdot 
\begin{pmatrix} 
0 \\ 1 \\ -\frac{1+p}{p} \cdot (b_{1} - b_{3}) \\ 1 
\end{pmatrix} + 
\frac{b_{1} - b_{3}}{2} \cdot 
\begin{pmatrix} 
0 \\ 0 \\ -\frac{1}{p} \cdot b_{1} - \frac{1}{2} (b_{1} - b_{3}) \\ 1 
\end{pmatrix}, 
\end{align*}
the cocycle 
$\beta_{ {\rm Sym}^{2}(\rho_{1}) }$
obtained from the second symmetric power representation 
${\rm Sym}^2(\rho_{1}) \: \widetilde{G}_{S}(\Q) \to \GL(3; \overline{\Q_{p}} )$ of 
the representation $\rho_{1}$ for the case of \cref{eg:cohAbel1+pGL3} 
is contained 
in the trivial cohomology class in 
$H^{1}(\widetilde{G}_{S}(\Q), ( \overline{\Q_{p}}^{\oplus 2} )_{\rho_{1}} )$, 
namely 
$$ 
[ \beta_{ {\rm Sym}^{2}(\rho_{1}) } ] = 
0 \in H^{1}(\widetilde{G}_{S}(\Q), ( \overline{\Q_{p}}^{\oplus 2} )_{\rho_{1}} ),
$$
and the representation $\rho_{1}$ for the case of 
\cref{eg:cohAbel1GL3} and \cref{eg:cohNonAbelGL3} 
is contained 
in the non-trivial cohomology class in 
$H^{1}(\widetilde{G}_{S}(\Q), ( \overline{\Q_{p}}^{\oplus 2} )_{\rho_{1}} )$, namely 
$$ 
[ \beta_{ {\rm Sym}^{2}(\rho_{1}) } ] \neq 
0 \in H^{1}(\widetilde{G}_{S}(\Q), ( \overline{\Q_{p}}^{\oplus 2} )_{\rho_{1}} ). 
$$
}
\end{rem}

\begin{eg} \label{eg:cohAbel1/1+pGL3} 
{\rm
For the same Galois group $\widetilde{G}_{S}(\Q)$ and 
the representation $\rho_{1} \: \widetilde{G}_{S}(\Q) \to \GL(2; \overline{\Q_{p}} )$
in \cref{eg:GL2}, 
assume $b_1 = b_3 = b$. 
Namely, we consider extensions of the homomorphism 
$\varphi_{2} \: F_{3} \to \GL(3; \overline{ \Q_{p}} )$. 
Consider 
the homomorphism  
$\widetilde{\varphi_{1}} \: F_{3} \to \GL(2; \overline{\Q_p} \otimes_{\Z_p} \Z_p[[\Gamma]] )$
defined by 
$$
\widetilde{\varphi_{1}} (g_{1}) =
\widetilde{\varphi_{1}} (g_{2}) =
\widetilde{\varphi_{1}} (g_{3}) = 
\begin{pmatrix} 
(1+p) \cdot \gamma & b \cdot \gamma \\ 
0 & \gamma 
\end{pmatrix}. 
$$
Then the matrix $Q_{ \widetilde{ \varphi_{1} } }$ is given by 
$$
Q_{\widetilde{ \varphi_{1} } } = 
\begin{pmatrix} 
 -p^{2} & 0 & p^{2} & 0 & 0 & 0 \\
 0 & -p^{2} & 0 & p^{2} & 0 & 0 \\
 p & 0 & -p & 0 & 0 & 0 \\
 0 & p & 0 & -p & 0 & 0 \\
 (1 + p) \gamma - 1 & b \gamma & - (1 + p) \gamma + 1 & -b \gamma & 0 & 0 \\
 0 & \gamma-1 & 0 & -\gamma + 1 & 0 & 0 \\
 - (1 + p) \gamma + (1 + 2p) & 
 -b \gamma & -p & 0 & (1 + p) (\gamma-1) & 
 b \gamma \\
 0 & - \gamma + (1 + 2p) & 0 & -p & 0 & \gamma - (1+p) \\
\end{pmatrix}, 
$$
and the GCD $\Delta_{2}(A_{\widetilde{ \varphi_{1} } })$ 
of all $(6-2)$-minors of $Q_{\widetilde{ \varphi_{1} } }$ is given by 
$$
\Delta_{2}(A_{\widetilde{ \varphi_{1} } }) \doteq 
(\gamma - 1) \{ \gamma - (1 + p) \} \in \Z_p[[\Gamma]]. 
$$

Let $a_{2} \neq 0, \ 1, \ 1+p$. 
Then 
$\mathbf{b} = {^t (b_{11}, b_{12}, b_{21}, b_{22}, b_{31}, b_{32} )}$ is a solution of 
$Q_{\widetilde{ \varphi_{1} } } |_{\gamma = a_{2}} \cdot \mathbf{b} = \boldsymbol{0}$ if and only if  
$b_{11} = b_{21} = b_{31}$, and 
$b_{12} = b_{22} = b_{32}$, namely, 
$$
\mathbf{b} = 
b_{11} \cdot 
\begin{pmatrix} 
1 \\ 0 \\ 1 \\ 0 \\ 1 \\ 0 
\end{pmatrix} + 
b_{12} \cdot 
\begin{pmatrix} 
0 \\ 1 \\ 0 \\ 1 \\ 0 \\ 1 
\end{pmatrix}. 
$$
Therefore, the homomorphism 
$\varphi_{2} \: F_{3} \to \GL(3; \overline{ \Q_{p}} )$ 
defined by 
$$
\varphi_{2}(g_{1}) = 
\varphi_{2}(g_{2}) =
\varphi_{2}(g_{3}) = 
\begin{pmatrix} 
a_{2} \cdot (1+p) & a_{2} \cdot b & b_{11} \\ 
0 & a_{2} & b_{12} \\ 
0 & 0 & 1 
\end{pmatrix}, 
$$
can be extended to 
a representation 
$\rho_{2} \: \widetilde{G}_{S}(\Q) \to \GL(3; \overline{\Q_{p}} )$. 

Let $a_{2} = \frac{1}{1+p}$, which is not a zero of 
$\Delta_{2}(A_{ \widetilde{ \varphi_{1} } })$, 
but the inverse is a zero of $\Delta_{2}(A_{ \widetilde{ \varphi_{1} } })$. 
Let $b = 0$. 
Consider the representation 
$\rho_{2} \: \widetilde{G}_{S}(\Q) \to \GL(3; \overline{\Q_{p}} )$,  
which is defined by 
$$
\rho_{2}(g_{1}) =
\rho_{2}(g_{2}) =
\rho_{2}(g_{3}) = 
\begin{pmatrix} 
1 & 0 & b_{11} \\ 
0 & \frac{1}{1+p} & b_{12} \\ 
0 & 0 & 1 
\end{pmatrix}. 
$$
Let $\beta_{2, (b_{11}, b_{12})} \: \widetilde{G}_{S}(\Q) \to \overline{\Q_{p}}^{\oplus 2}$ 
be a 1-cocycle with coefficient $\varphi_{2}$
such that 
$\beta_{2, (b_{11}, b_{12})} (g_{1}) = 
\beta_{2, (b_{11}, b_{12})} (g_{2}) = 
\beta_{2, (b_{11}, b_{12})} (g_{3}) = {}^{t}(b_{11}, b_{12})$. 

If $b_{11} \neq 0$ and $b_{12} = 0$, 
then $\beta_{2, (b_{11}, 0)}$ is not a 1-coboundary 
since we cannot take 
$v = {}^{t}(v_{1}, v_{2}) \in \overline{\Q_{p}}^{\oplus 2}$ which satisfies \eqref{eq:cobdry}, 
namely 
$$
\displaystyle
\begin{pmatrix} 
(1 - 1) \cdot v_{1} + 0 \cdot v_{2} \\ 
\left( \frac{1}{1+p} - 1 \right) \cdot v_{2} 
\end{pmatrix}
= 
\begin{pmatrix} 
b_{11} \\ 
0  
\end{pmatrix}. 
$$

If $b_{11} = 0$ and $b_{12} \neq 0$, 
then $\beta_{2, (0, b_{12})}$ is a 1-coboundary 
since we can take 
$v = {}^{t}(v_{1}, v_{2}) \in \overline{\Q_{p}}^{\oplus 2}$ as 
$v = 
\displaystyle 
\begin{pmatrix} 
1 \\ 
-\frac{1+p}{p} \cdot b_{12} 
\end{pmatrix}
$ 
which satisfies \eqref{eq:cobdry}, 
namely 
$$
\displaystyle
\begin{pmatrix} 
(1 - 1) \cdot v_{1} + 0 \cdot v_{2} \\ 
\left( \frac{1}{1+p} - 1 \right) \cdot v_{2} 
\end{pmatrix}
= 
\begin{pmatrix} 
0 \\ 
b_{12} 
\end{pmatrix}. 
$$

Therefore, we have 
\begin{align*}
\dim_{\overline{\Q_{p}}} Z^{1}(\widetilde{G}_{S}(\Q), ( \overline{\Q_{p}}^{\oplus 2} )_{\rho_{1}} ) ) &= 2, \\
\dim_{\overline{\Q_{p}}} B^{1}(\widetilde{G}_{S}(\Q), ( \overline{\Q_{p}}^{\oplus 2} )_{\rho_{1}} ) ) &= 1, \\
\dim_{\overline{\Q_{p}}} H^{1}(\widetilde{G}_{S}(\Q), ( \overline{\Q_{p}}^{\oplus 2} )_{\rho_{1}} ) ) &= 1. 
\end{align*}
}
\end{eg}


\noindent
{\it Acknowledgements.}
The authors would like to thank 
Takahiro Kitayama, 
Masanori Morishita,
Yoshikazu Yamaguchi 
and the anonymous referee 
for helpful communications. 
This work is partially supported by JSPS KAKENHI Grant Numbers JP22K03268, JP22H01117, and JP21K03240. 


\printbibliography


\ 

\noindent
Yasushi Mizusawa mizusawa@rikkyo.ac.jp \\ 
Department of Mathematics, Rikkyo University, 3-34-1 Nishi-Ikebukuro, Toshima-ku, 171-8501 Tokyo, Japan

\

\noindent
Ryoto Tange rtange.math@gmail.com \\ 
Department of Mathematics, School of Education, Waseda University, 1-104, Totsuka-cho, Shinjuku-ku, Tokyo, 169-8050, Japan

\

\noindent
Yuji Terashima yujiterashima@tohoku.ac.jp \\ 
Graduate School of Science, Tohoku University, 6-3, Aoba, Aramaki-aza, Aoba-ku, Sendai, 980-8578, Japan

\end{document}